\documentclass[11pt]{article}
\usepackage{amssymb,psfig,latexsym}

\newtheorem{theorem}{Theorem}[section]
\newtheorem{lemma}{Lemma}[section]
\newtheorem{exam}{Example}[section]

\newcommand{\bsq}{\vrule height .9ex width .8ex depth -.1ex}

\newcommand{\hsp}{\hspace{\parindent}}

\newcommand{\RR}{{\Bbb R}}
\newcommand{\QQ}{{\Bbb Q}}

\newcommand{\ZZ}{{\Bbb Z}}

\newcommand{\sC}{{\cal C}}

\newcommand{\sO}{{\cal O}}

\newcommand{\bv}{{\bf v}}
\newcommand{\bw}{{\bf w}}

\newcommand{\beql}[1]{\begin{equation}\label{#1}}
\newcommand{\eqn}[1]{(\ref{#1})}
\newcommand{\eeq}{\end{equation}}

\catcode`\@=11
\renewcommand{\section}{
        \setcounter{equation}{0}
        \@startsection {section}{1}{\z@}{-3.5ex plus -1ex minus
        -.2ex}{2.3ex plus .2ex}{\large\bf}
        }
\catcode`@=12
\makeatletter
\def\eqalignno#1{\displ@y \ta {\bf s} kip\@centering
  \halign to\displaywidth{\hfil$\@lign\displaystyle{##}$\ta {\bf s} kip\z@skip
    & $\@lign\displaystyle{{}##}$\hfil\ta {\bf s} kip\@centering
    & \llap{$\@lign##$}\ta {\bf s} kip\z@skip\crcr
    #1\crcr}}
\makeatother

\makeatletter
\def\@sect#1#2#3#4#5#6[#7]#8{\ifnum #2>\c@secnumdepth
     \def\@svsec{}\else 
     \refstepcounter{#1}\edef\@svsec{\csname the#1\endcsname.\hskip .75em }\fi
     \@tempskipa #5\relax
      \ifdim \@tempskipa>\z@ 
        \begingroup #6\relax
          \@hangfrom{\hskip #3\relax\@svsec}{\interlinepenalty \@M #8\par}%
        \endgroup
       \csname #1mark\endcsname{#7}\addcontentsline
         {toc}{#1}{\ifnum #2>\c@secnumdepth \else
                      \protect\numberline{\csname the#1\endcsname}\fi
                    #7}\else
        \def\@svsechd{#6\hskip #3\@svsec #8\csname #1mark\endcsname
                      {#7}\addcontentsline
                           {toc}{#1}{\ifnum #2>\c@secnumdepth \else
                             \protect\numberline{\csname the#1\endcsname}\fi
                       #7}}\fi
     \@xsect{#5}}
\def\@theorem#1#2{\it \trivlist \item[\hskip \labelsep{\bf #1\ #2.}]}
\makeatother

\begin{document}
\begin{center}
{\large {\bf Dynamics of  a Family of Piecewise-Linear 
Area-Preserving Plane Maps}} \\
{\large {\bf II. Invariant Circles}} \\ \bigskip
{\large {\em Jeffrey C. Lagarias}} \\ \smallskip
Department of Mathematics \\
University of Michigan  \\
Ann Arbor, MI 48109-1043 \\
{\tt email:}  {\tt lagarias@umich.edu}\bigskip \\

{\large {\em Eric Rains}} \\ \smallskip
Department of Mathematics \\
University of California-Davis \\
Davis, CA 95616-8633  \\
{\tt email:}  {\tt rains@math.ucdavis.edu}\bigskip \\
\vspace*{2\baselineskip}
(July 1, 2005 version) \\

\vspace*{1.5\baselineskip}
{\bf ABSTRACT}
\end{center}

This paper studies the behavior under iteration of 
the  maps $T_{ab}(x, y) = (F_{ab}(x) - y, x)$ of the
plane $\RR^2$, in which $F_{ab}(x) = ax$ if $ x \ge 0$ and $bx$ if $x < 0$.
The orbits under
iteration correspond to solutions of the difference equation
$x_{n+2}= 1/2(a-b)|x_{n+1}| + 1/2(a+b)x_{n+1} - x_n.$
This family of piecewise-linear
maps of the plane has the  parameter space $(a, b) \in \RR^2$. 
These maps are area-preserving homeomorphisms of $\RR^2$ that map
rays from the origin into rays from the origin. 
We show the existence of  special parameter values  where $T_{ab}$
has every nonzero orbit contained in an  invariant
circle with an irrational rotation number, with 
 invariant circles that are  piecewise unions of arcs of conic sections.
Numerical experiments indicate the possible existence of 
invariant circles for many other parameter values.

\vspace*{1.5\baselineskip}

\setlength{\baselineskip}{1.0\baselineskip}

%
%
%
%
%

\section{Introduction}
\hsp

We continue the study the behavior under iteration of the two parameter 
family of piecewise-linear
homeomorphisms of $\RR^2$ given by
\beql{eq101}
T_{ab} (x,y) =
\left\{
\begin{array}{ccc}
(ax-y,x) & \mbox{if} & x \ge 0 , \\
(bx-y,x) & \mbox{if} & x < 0 .
\end{array}
\right.
\eeq
The parameter space is $(a,b) \in \RR^2$.
This map can be written
\beql{eq102}
T_{ab} (x,y) = \left[
\begin{array}{cc}
F_{ab} (x) & -1 \\ 1 & 0
\end{array}
\right]
\left[
\begin{array}{c}
x\\y
\end{array}\right]
 ~,
\eeq
in which
\beql{eq103}
F_{ab} (x) = \left\{
\begin{array}{ccc}
a & \mbox{if} & x \ge 0 , \\
b & \mbox{if} & x < 0;
\end{array}
\right.
\eeq
we view elements of $(x,y)\in\RR^2$ as column vectors.
The formula \eqn{eq102} shows that $T_{ab} (x,y)$ is a homeomorphism, 
since
\begin{eqnarray}\label{eq104}
T_{ab}^{-1} (x,y) & = &
\left[ \begin{array}{cc}
F_{ab} (y) & -1 \\ 1 & 0
\end{array}\right]^{-1} 
\left[\begin{array}{c}x\\y\end{array}\right] 
\nonumber \\
& = & \left[ \begin{array}{cc}
0 & 1 \\ -1 & F_{ab} (y)
\end{array}\right]
\left[\begin{array}{c}x\\y\end{array}\right] 
.
\end{eqnarray}
It preserves the area form $d \omega = dx \wedge dy$, and 
it also maps rays from the origin into rays 
from the origin. The orbits of the iteration are $T_{ab}(x_{n+1}, x_{n})=
(x_{n+2}, x_{n+1})$, in which 
the quantities $ x_{n}$ satisfy the difference equation
\beql{eq101c}
x_{n+2}= \mu|x_{n+1}| + \nu x_{n+1} - x_n.
\eeq
in which $\mu= 1/2(a-b)$ and $\nu= 1/2(a+b)$. For certain
purposes the $(\mu, \nu)$ coordinate system is more convenient
than the $(a,b)$ coordinate system (see part III), but 
in this paper we shall exclusively use the $(a,b)$ coordinate
system.

The main result of this paper  is the characterization
of certain one-parameter families of $T_{ab}$
whose generic members exhibit quasiperiodic motion
under iteration, so have
bounded orbits   (Theorem ~\ref{th25}).
These families are those 
for which $(0,1)$ and 
$(0,-1)$ fall in the same orbit with a fixed
symbolic dynamics between their occurrences,
and which contain a parameter value
having an irrational rotation number. We show that
generic members of any such family have each orbit 
contained in an invariant circle which is a 
finite union of arcs of conic sections. We
rigorously establish for specific examples 
that the irrational rotation
number property holds (Examples \ref{ex51} and \ref{ex52}). 

Recall that $\Omega_{SB}$ is the set of all parameter
values for which the map $T_{ab}$ has at least one
nonzero  bounded orbit, $\Omega_{B}$ is the the set of
paramenter values having all orbits bounded
and  $\Omega_{R}$
is the set of parameter values for which
$T_{ab}$ is topologically conjugate to a rotation
of the plane, or equivalently, having an invariant circle.
(We regard a purely periodic map as having an invariant
circle, which may be constructed artificially from
a suitable union of orbits, as all orbits
are finite.)
Our results here have a bearing on 
$\Omega_{SB}$, because 
there  appears to be a close relation between 
the narrow set of parameter values having piecewise
conic invariant circles, which we show has
two-dimensional Lebesgue measure zero,   and the larger
set $\Omega_{SB}$. 
Let $\Omega_{P}$ denote the set of parameter
values $(a,b)$ for which $T_{ab}$ is periodic, and let
 $\Omega_{Q}$ ( Q = ``quadratic'') 
denote the set of parameter values
having a piecewise conic invariant circle as
prescribed by Theorem~\ref{th55}.
Then we have the inclusions
$$
\Omega_{P} \subset \Omega{P} \cup \Omega_{Q} \subseteq \Omega_{R}
=  \Omega_{B} \subseteq \Omega_{SB}.
$$
We formulate the following conjecture. \\

\noindent{\bf Conjecture A.}
{\em The set $\Omega_{P} \cup \Omega_{Q}$ of pairs
$(a,b)$ for which $T_{ab}$ is either periodic or has a 
piecewise conic invariant circle is dense in $\Omega_{SB}$.} \\

A closely related conjecture was made by Bedford, Bullett and
Rippon \cite{BBR95}, to the effect that the 
parameter set of purely periodic maps 
$\Omega_P$ is dense in $\Omega_{B}$. In part III we will show
that $\Omega_{SB}$ is a closed set, and it is natural to
expect that it is the closure of $\Omega_{B}$, though
we have not established this.  Combining this expectation with 
the conjecture of \cite{BBR95}, would result in
 the conjecture
that $\Omega_{P}$ is dense in $\Omega_{SB}$.

\noindent
In \S5 we show that  $\Omega_{P} \cup \Omega_{Q}$ has
Hausdorff dimension $1$, hence has two-dimensional Lebesgue measure zero. 
We also  present nonrigorous
numerical evidence in \S4 of invariant circles for
a variety of parameter values which appear to fall
outside the set  $\Omega_{P} \cup \Omega_{Q}$.

There has been previous work concerning  invariant circles
for these maps.
In 1986 M. Herman \cite[Theorem VIII.5.1]{He86} proved 
results implying that
any map $T_{ab}$  having an irrational rotation number $r(S_{ab})=r$
such that $r$ has bounded partial quotients in its continued fraction
expansion is  topologically conjugate to a
rotation of the plane, so has an invariant circle. 
Herman's results also supply some
justification for the numerical observations that
are noted in \S4; his set of allowed parameter values
presumably has Hausdorff dimension $2$, but with 
two-dimensional Lebesgue measure zero.
It undoubtably include  
parameter values outside $\Omega_{P} \cup \Omega_{Q}$.

In 1995 Beardon, Bullett and Rippon \cite[p.673]{BBR95} announced 
without proof as work to come
 results related to those appearing here. 
They asserted that maps $T_{ab}$ mapping $(0, 1)$  to $(0,-1)$
after a finite number of steps with a fixed symbolic itinerary 
should give an algebraic curve of parameter values on which
there is an open subset of values 
where the map is topologically conjugate to 
a rotation of the plane, compare our Theorem~\ref{th25}.
Example~\ref{ex53} below shows that their assertion
requires additional  hypotheses 
(such as  nonconstancy of  the rotation number
along the curve) 
for the stated conclusion 
to hold. It seems clear they had knowledge 
of various ideas
presented here; however they made no
comment regarding the nature of the 
invariant circles.

\noindent \paragraph{Notation.} We write $\bv=(\bv_x, \bv_y) \in \RR^2$,
to be viewed as a column vector.
 An interval $[\bv_1, \bv_2)$ of the 
unit circle, or corresponding sector $\RR^{+}[\bv_1, \bv_2)$
of the plane $\RR^2$, is the one specified by going counterclockwise from
$\bv_1$ to $\bv_2$.

\noindent \paragraph{Acknowledgments.} We did most of 
the work reported in this paper 
while employed at AT\&T Bell Labs; most results of this
paper were obtained during the summer of 1993.
We thank M. Kontsevich for
bringing the work of Beardon, Bullett and Rippon \cite{BBR95}
to our attention.

%
%
%
%
%

\section{Summary of Results}
\hsp
It is convenient to represent the action of $T_{ab}$, 
acting on row vectors $\bv_n = (x_{n+1}, x_n )$ as
\beql{eq202}
T_n (\bv_0) = (x_{n+1}, x_n)
 = M_n (\bv_0) (x_1,  x_0)
\eeq
in which
\beql{eq203}
M_n (\bv_0) = \prod_{i=1}^n \left[
\begin{array}{cc}
F_{ab} (x_i) & -1 \\ 1 & 0
\end{array}
\right] :=
\left[
\begin{array}{cc}
F_{ab} (x_n) & -1 \\ 1 & 0
\end{array} 
\right] \cdots
\left[
\begin{array}{cc}
F_{ab} (x_2 ) & -1 \\ 1 & 0 
\end{array}
\right] 
\left[
\begin{array}{cc}
F_{ab} (x_1 ) & -1 \\
1 & 0
\end{array}
\right]\,.
\eeq

Conjugation by the involution $J_0: (x,y) \to (-x, -y)$ gives
\beql{eq204}
T_{ba} = J_0^{-1} \circ T_{ab} \circ J_0 ~.
\eeq
Thus, in studying dynamics, without loss of generality we can
restrict to the closed half-space
$\{(a,b) : a \ge b \}~$ of the $(a,b)$ parameter space.

The associated rotation map $S_{ab}: S^1 \to S^1$ is given
by 
$$ S_{ab}(e^{i\theta}):= 
\frac{T_{ab}(e^{i\theta})}{|T_{ab}(e^{i\theta})|}.
$$
It has a well-defined rotation number $r(S_{ab})$, which
is counterclockwise rotation, and was shown in part I to
always lie in the closed interval $[0, 1/2].$

In part I we classified the dynamics of $T_{ab}$ 
 in the case where the
rotation number $r(S_{ab})$ is rational, as
follows \cite[Theorem 2.4]{LR02aa}.

\begin{theorem}\label{th24}
If the rotation number
$r(S_{ab} )$ is rational, then $S_{ab}$ has a periodic orbit,
and one of the following three possibilities occurs.

(i) $S_{ab}$ has exactly one periodic orbit.
Then $T_{ab}$ has exactly one periodic orbit (up to scaling) 
and all other orbits diverge in modulus to $+ \infty$ as $n \to \pm \infty$.

(ii) $S_{ab}$ 
has exactly two periodic orbits.
Then $T_{ab}$ has no periodic orbits.
All orbits of $T_{ab}$ diverge in modulus
to $+\infty$ as $n \to \pm \infty$, with  
the exception of orbits lying over the two periodic orbits of $S_{ab}$.
These exceptional orbits have  modulus diverging  to $+ \infty$ 
in one direction and to $0$ in the other direction, with forward
divergence for one, and backward divergence for the other.

(iii) $S_{ab}$ has at least three periodic orbits.
Then $T_{ab}$ is of finite order, i.e. $T_{ab}^{(k)} =I$ for some $k \ge 1$,
and all its orbits are periodic.
\end{theorem}

In this paper we are mainly concerned with parameter values
where the rotation number is irrational, but we will need
the result above.
The main result of this
paper is that there exist parameter values $(a,b)$
with irrational rotation number
for which $T_{ab}$ has  invariant circles 
with a striking structure.

\begin{theorem}\label{th25}
Suppose that the rotation number $r(S_{ab})$ is irrational,
and that the $S_{ab}$ orbit of $(0,1)$ contains $(0,-1)$.
Then the following hold.

(1) The $T_{ab}$ orbit of $(0,1)$ contains $(0,-1)$. 

(2) The closure of every (nonzero) orbit of $T_{ab}$ 
is an invariant circle, 
which is a piecewise union of arcs of  conic sections.
The conic sections occurring in such an invariant circle
are all of the same type, either
ellipses, hyperbolas or straight lines. 
\end{theorem}

This result follows from Theorem~\ref{th55},
which gives more detailed information about the
 the number of conic pieces in such an invariant circle. \\

In \S4 we give examples exhibiting 
 invariant circles of all the types allowed by
Theorem~\ref{th25}. 
An example is $a={2}^{1/4}$, $b = - {2}^{1/4}$,
which has irrational rotation number and
whose  invariant circles are the 
union of eight segments of ellipses.
We actually exhibit one-parameter families
whose generic members have the required properties.
Here ``generic'' means irrational rotation number;
these one-parameter families include values 
$T_{ab}$ with rational rotation
number, where $T_{ab}$ is a  periodic map.
In the Appendix we  verify the conditions
of Theorem ~\ref{th55} hold in these examples, and
exhibit  some specific  parameter values where
the rotation number is irrational. 
In some cases we are able to show
the rotation number is transcendental,
using Baker's method in transcendental
number theory.
In \S4 we also give a variety of computer plots
for parameter values not covered by Theorem~\ref{th25}
which numerically appear to  give invariant circles.
We have no rigorous proof of the existence of
non piecewise-conic invariant circles, however. 

In \S5 we show that the set
$\Omega_{Q}$ has Hausdorff dimension $1$
and discuss supporting evidence for Conjecture A.
We conclude with a discussion of the possiblity
of $T_{ab}$ with
associated circle map  having irrational rotation number, but
with $T_{ab}$ not being conjugate to a rotation of the plane.

%
%
%
%
%
\section{Irrational Rotation Number}
\hsp

In this paper we  study maps in the irrational rotation case 
by considering their first-return maps to
suitable sectors of the plane.
%
%
\begin{theorem}\label{th53}
Suppose that the rotation number $r(S_{ab} )$ is irrational.

(1) For any half-open sector $J= \RR^{+}[\bv , \bv' )$ 
the first return map $T_J^{(1)}: J \to J$ of
$T_{ab}$ to $J$
 is piecewise linear with at most five pieces.

(2) Let $m_+$, $m_-$, $n_+$, $n_-$ be nonnegative integers, 
and set 
$\sO := \sO_1 \cup \sO_2 \cup \sO_3 \cup \sO_4$, where:
\begin{eqnarray*}
\sO_1 & = & \{ 
T_{ab}^{(1+i)} (0,1) : 0 \le i < m_+ \} \\
\sO_2 & = & \{ T_{ab}^{(-i)} (0,1) : 0 \le i < m_- \} \\
\sO_3 & = & \{ T_{ab}^{(1+i)} (0,-1): 0 \le i < n_+ \} \\
\sO_4 & = & \{T_{ab}^{(-i)} (0,-1) : 0 \le i < n_- \} \,.
\end{eqnarray*}
Let $J = \RR^{+}[\bv, \bv' )$ be any half-open sector determined by two
elements $\bv , \bv' \in \sO$  such that $J$ contains no points of
$\sO$ in its interior.
Then the first return map $T_J : J \to J$ is piecewise linear 
with at most three pieces.
\end{theorem}

\noindent\paragraph{Proof.}
(1) Note first that for any half-open sector $J$, the first
return map for $T_{ab}$ is well-defined,
since  the assumption that $S_{ab}$ has
irrational rotation number implies 
that the orbit of any point under $S_{ab}$ is dense.

Let $\bv_1$ be the first preimage $T_{ab}^{(-i)}(0,1)$, $i \ge 0$,    
of $(0, 1)$  that falls in $J$;
let $\bv_2$ be the first  preimage of $(0,-1)$ with the same property. Let
$\bv_3$ be the first point 
of the form $T_{ab}^{(-i)}(\bv)$, $i>0$ that falls in $J$,
 and let
$\bv_4$ be the first point of the form $T_{ab}^{(-i)}(\bv')$, $i\ge 0$
that falls in $J$.  The $\bv_i$ split
$J$ into at most five intervals. We claim that the first return map
for $T_{ab}$ on the corresponding sector is linear on each corresponding
subsector.

Consider first the (at most) three half-open subsectors into which
$J$ is subdivided by 
$\bv_3$ and $\bv_4$.
The number of steps until the first return will be constant 
on each subsector,
since between any two points that take a different number of steps to
return, there must be a point that hits the boundary of $J$ at or before the
first return.  If $\bv_5$ were such a point other than $\bv_3$ or $\bv_4$, 
then the
iterates of $\bv_5$ would hit the ray determined by
$\bv_3$ or $\bv_4$ before hitting the ray determined by
$\bv$ or $\bv'$, which gives
a contradiction, since hitting $\bv_3$ or $\bv_4$ constitutes 
a return to $J$.

Now, we restrict our attention to one of these three subsectors,  
on which the number
of iterations  until return to $J$ is constant. 
For any two points in the subsector whose iterates
until first return to $J$ are always on the same side of the $y$-axis, the
first return map on $J$ is linear on the sector they induce.  
Thus the only points where linearity can fail are points which
hit the $y$-axis on iteration before their
first return, namely  $\bv_1$ and $\bv_2$. (Any other such
point on $J$ would hit the ray determined by
$\bv_1$ or $\bv_2$ before it hit the $y$-axis, a contradiction.)
Thus, the subsector  between each two neighboring such $\bv_i$ 
is a subsector on which
the first return map on $J$ is linear, and  (1) follows.

(2) First note that $T_{ab}(0,1) = (-1,0)$; 
similarly $T_{ab}(0, -1) = (1,0)$.
Suppose now that $\bv \in \sO_1$. by definition $\bv_3$ is the first
strict preimage
of $\bv$ that falls in $J$. This preimage cannot belong to  $\sO_1$ 
because by hypothesis the only point of $\sO$ in the half-open
interval $J$ is the endpoint
$\bv$. Thus this preimage is a strict preimage of $(-1,0)$, 
hence is $(0,1)$ or a preimage, so we
conclude that $\bv_3 = \bv_1$. 

We deduce the following facts by similar reasoning:
(1) If $\bv \in \sO_2$ then $\bv_1 = \bv$;
(2) If $\bv \in \sO_3$ then $\bv_3 = \bv_2$;
(3) If $ \bv \in \sO_4$ then $\bv_2 = \bv$.

The preimage $\bv_4$ of $\bv'$ must also be
equal to one of these other values.  
We distinguish three cases: (a) $\bv$ and $\bv'$ are in the same orbit,
with $\bv$ coming first, (b) $\bv$ and $\bv'$ are in the same orbit, with
$\bv'$ coming first, (c) $\bv$ and $\bv'$ are in different orbits.
In case (a), $\bv_4=\bv$, while in case (b), $\bv_4=\bv_3$.
In case (c) one deduces 
that if $\bv' \in \sO_1 \cup \sO_2$ then $\bv_4 = \bv_1$, 
while if $\bv' \in \sO_3 \cup \sO_4$ then $\bv_4 = \bv_2$.
In all cases the five points $\bv$, $\bv_1$, $\bv_2$, $\bv_3$, $\bv_4$
take at most three different values, and therefore $T_J$ 
has at most three linear pieces on $\RR^{+} [J]$.~~~$\bsq$ \\

We next  consider the special case where the 
 $T_{ab}$-orbit of $(0,1)$ contains
$(0,-1)$. To analyze this case
 we  recall the following result from part I 
\cite[Theorem 3.4]{LR02aa}.
%
%

\begin{theorem}~\label{th34}
Let $n\ge 1$ and suppose that
$S_{ab}^{(n)}(0,1) = (0, \pm 1),$
or $S_{ab}^{(n)}(0,-1) = (0, \pm 1),$
Then one of the following two relations holds:
\beql{eq400a}
T^{(n)}_{ab}(0,1)=(0,\lambda),
\eeq
\beql{eq400b}
T^{(n)}_{ab}(0,-1)=(0, - \lambda^{-1}).
\eeq
where $\lambda$ is a nonzero real number.

(i) If $\lambda>0$, then 
both relations above hold. In addition, 
\beql{307a}
T^{(n)}_{ab}(-1,0) =  (-\lambda,0) ~~~~\mbox{and}~~~~
T^{(n)}_{ab}(1,0) =  (\lambda^{-1},0).
\eeq
The rotation number $r(S_{ab})$ is rational.

(ii) If $\lambda < 0 $, then necessarily $\lambda=-1$.
In the first case
\beql{308a}
T^{(n)}_{ab}(0, 1)   =  (0,-1) ~~~~\mbox{and}~~~~
T^{(n)}_{ab}(- 1,0)  =  (1,0),
\eeq
while in the second case,
\beql{309a}
T^{(n)}_{ab}(0,-1)  =  (0,1) ~~~~\mbox{and}~~~~
T^{(n)}_{ab}(1,0)   =  (-1,0).
\eeq
The rotation number $r(S_{ab})$ can be irrational or rational.
\end{theorem}

We now show that the $T_{ab}$-orbit of $(0,1)$ contains
$(0,-1)$, then one can obtain sectors 
on which the first return map is piecewise linear
with at most two pieces, 
and a surprising thing happens:
%
%
\begin{theorem}\label{th54}
Suppose that $S_{ab}$ has irrational rotation number, and that $(0,1)$
and $(0, -1)$ are in the same $T_{ab}$-orbit. 
Let $\sO = \{ T_{ab}^{(j)}(\bv_0): 1 \le j \le m \}$ be any finite segment
of that orbit which includes all of the points strictly between
$(0,1)$ and $(0, -1)$, as well as whichever of $(0,1)$ and $(0,-1)$ comes
last in the orbit.  If $J = \RR^{+}[ \bv, \bv' )$ is a half-open sector
determined by two members of $\sO$ containing no elements of $\sO$ in its
interior, then the first return map $T_J$ on $J$ is 
piecewise linear with exactly
two pieces, and the linear transformations corresponding to the two pieces
commute.
\end{theorem}

\paragraph{Proof.}
Let $J= \RR^{+}[\bv, \bv')$
be a sector determined by two members of $\sO$, such that no element
of $\sO$ lies in its interior. We show that the first return map $T_J$
of $T_{ab}$ on the sector $J$ 
 is piecewise linear with exactly two pieces.
Let $T_{ab}^{(n)}(0,1)= (0, -1)$.
We first treat the case $n > 0$, when $(0,1)$ comes earlier
in the orbit than $(0,-1)$. In the notation
of Theorem ~\ref{th53}, we take $\sO_i$ as large as possible inside $\sO$,
and obtain the following relations:
$\sO_2\subset\sO_4$, $\sO_3\subset\sO_1$, 
and $\sO_1\cup\sO_2=\sO_3\cup\sO_4$.
Theorem~\ref{th53}(ii) states that the first return map $T_J$ is piecewise
linear with at most three pieces. We show
that the relations between the sets $\sO_i$
cause two of the pieces to become identified.  Indeed, from the proof of
Theorem~\ref{th53}(ii), we have the following three possibilities:

(a) $\bv\in \sO_1\cap\sO_4$: $\bv_3=\bv_1$, $\bv_2=\bv$.

(b) $\bv \in \sO_2\cap\sO_4$: $\bv_2=\bv_1=\bv$.

(c) $\bv \in \sO_3\cap\sO_1$: $\bv_3=\bv_2=\bv_1$.

\noindent In case (a), $\bv_4=\bv$, while in case (b), $\bv_4=\bv_3$; 
case (c) is
excluded, since $\bv$ and $\bv'$ are in the same orbit.  Thus the five
points $\bv$, $\bv_1$, $\bv_2$, $\bv_3$, $\bv_4$ take on only two values
($\bv$ and $\bv_3$), which divides the sector $J$ into two pieces. We
cannot have a further collapse to one piece because the equality
$\bv=\bv_3$ would make $\bv$ a periodic point of $S_{\mu\nu}$, whence
$S_{\mu\nu}$ has a rational rotation number, contradicting the irrational
rotation hypothesis.  We conclude that $T_J$ is piecewise linear on the
sector $J$ with exactly two pieces.

In the other case $n < 0$ we have
$\sO_4\subset\sO_2$, $\sO_1\subset\sO_3$, 
and $\sO_1\cup\sO_2=\sO_3\cup\sO_4$.
This yields the three possibilities:

(a) $\bv\in \sO_2\cap\sO_3$: $\bv_3=\bv_2$, $\bv_1=\bv$.

(b) $\bv\in \sO_2\cap\sO_4$: $\bv_2=\bv_1=\bv$.

(c) $\bv\in \sO_3\cap\sO_1$: $\bv_3=\bv_2=\bv_1$.

\noindent By similar arguments to the case $n>0$, we conclude that
 $T_J$ is piecewise linear with exactly two pieces.

Now let $J=\RR^{+}[\bv,\bv')$ be a sector of the
type above. The first return map  $T_J$
is piecewise linear on the sector with
exactly two pieces, and we let 
 $\bv^\ast$ denote the breakpoint in the interior
of $J$ (so $\bv^\ast=\bv_3$.)
Let $M_1$ be the linear transformation corresponding to 
$J_1=\RR^{+}[\bv,\bv^\ast)$,
and let $M_2$ be the linear transformation corresponding 
to $J_2= \RR^{+}[\bv^\ast,\bv')$, so that
$$
T_J(\bw) = \left\{ { {M_1\bw  ~~\mbox{if}~~\bw \in \RR^+[\bv, \bv^\ast)} 
\atop {M_2\bw  ~~\mbox{if}~~\bw \in \RR^{+}[\bv^\ast, \bv')} }\right..
$$ 
Note that $M_2 \bv^\ast =\bv$, and, 
by continuity, $M_1\bv^{\ast}=\bv'$.  Now, consider
the action of $T_J^{(2)}$.  In some neighborhood of $\bv^\ast$,
 $T_J^{(2)}$ acts
as $M_2M_1$ on points clockwise of $\bv^\ast$ and as $M_1M_2$ on points
counterclockwise of $\bv^\ast$.  If we can show that 
$T_J^{(2)}$ is linear on
this neighborhood, then $M_1M_2=M_2M_1$ and we are done.

     To compute $T_J^{(2)}$, one multiplies a sequence of matrices of the
form $\left[{c\atop 1}{-1\atop 0}\right]$, with $c=a$ or $c=b$ at step $i$
if the $i$th $x$-coordinate is positive or negative, respectively.  Now,
as we pass through $\bv^\ast$, the sequence will change only at the 
two steps $j=j_1, j_2$ where
$T^{(j)}_{ab}(\bv^\ast)=(0,\pm 1)$.  Furthermore, away from 
$\bv^\ast$ in the
neighborhood, if $c=a$ at one of those two steps, we must have $c=b$ at the
other.  Thus we find that (say)
$$M_1M_2 = M_A \left[{b\atop 1}{-1\atop 0}\right]
M\left[{a\atop 1}{-1\atop 0}\right] M_B$$
and
$$M_2M_1 = M_A \left[{a\atop 1}{-1\atop 0}\right]
M\left[{b\atop 1}{-1\atop 0}\right] M_B,$$
for certain matrices $M_A, M_B$ and $M$,  in which
$\tilde{M}:=M\left[{*\atop 1}{-1\atop 0}\right]$
is a product of matrices corresponding to the sequence
of iterations between
the two places where $T_{ab}^{(j)}(\bv^\ast)=(0,\pm 1)$.
Note that the value of $*$ ($a$ or $b$) has no effect on the first
step of the iteration of $T_{ab}$.
The  matrix equation $M_1M_2=M_2M_1$ boils down to showing that
$$
\left[{a\atop 1}{-1\atop 0}\right]M\left[{b\atop 1}{-1\atop 0}\right]=
\left[{b\atop 1}{-1\atop 0}\right]M\left[{a\atop 1}{-1\atop 0}\right].
$$

We now recall from Theorem~\ref{th34}(ii),
the key fact that  if
 $S_{ab}^{(n)}(0,1)=(0,-1)$ for some (positive or negative)
integer $n$, then
$$T_{ab}^{(n)}(0,1)=(0,-1).$$
If $n>0$, this gives 
$\tilde{M}(0,1) = M(-1,0) =(0,-1)$; if $n<0$ 
we obtain instead
$ \tilde{M}(0, -1) = M(1,0)= (0, 1)$.
In either case, we can conclude that
$$
M=\left[{0 \atop 1}{-1\atop d}\right],
$$
for some value $d$, and  it follows that
$$
\left[{a\atop 1}{-1\atop 0}\right]M\left[{b\atop 1}{-1\atop 0}\right]=
\left[{-a-b-d\atop -1}~~{1\atop 0}\right],
$$
This is symmetric in $a$ and $b$, hence $M_1M_2=M_2M_1$,
and the theorem is proved.~~~$\bsq$ \\

In preparation for the next result, we recall the following well-known
fact, see Arnold and Avez~\cite[Appendix 27]{AA68}.      .
%
%
\begin{lemma}\label{le51}
Let 
$M = \left[{a \atop c}{b \atop d}\right]   \in SL(2, \RR)$, 
and suppose that $M \neq \pm I$. Then the map 
 $(x,y) \to M(x, y)$
leaves invariant the quadratic form
$$Q(x,y) = cx^2 + (d-a)xy -by^2,$$
and the only quadratic forms it leaves invariant are scalar
multiples of $Q(x,y)$. The level sets
$  cx^2 + (d-a)xy -by^2 = \lambda $
for real $\lambda$ are invariant sets, and fill the plane.
Furthermore:

(1) If $|Tr(M)| < 2$, then 
the nonempty level sets  are ellipses, or else a single point.

(2) If $|Tr(M)| > 2 $, 
then the level sets are  hyperbolas, except for 
$\lambda=0$, where they are their asymptotes,
consisting of two straight lines through the origin. 

(3) If $|Tr(M)| = 2 $, then 
the  nonempty level sets are either two parallel lines
whose vector sum is a fixed parallel line through the origin;
or, for $\lambda=0$, this fixed line.

Any matrix $M'$ commuting with $M$ also preserves the quadratic form $Q$.
\end{lemma}

\paragraph{Proof.}
The condition for invariance of the quadratic form
$Q(x, y) = Ax^2 + 2Bxy +Cy^2$ is that $M^T QM = Q,$ where 
$Q =\left[{A \atop B}{B\atop C}\right]$; an
equivalent matrix condition is that $M^T Q= Q M^{-1}$.
Letting
$M=\left[{a \atop c}{b\atop d}\right]\in SL_2(\RR),$
the latter condition yields the linear equations
$(a-d)A = -2cB $,
$ cC = -bA $, and 
$(a-d)C = 2bB.$ At least one of $b, c, a-d$ is nonzero, otherwise
$M = \pm I$ which we have excluded. Then this linear system has rank 
at least two, so has at most
a one-parameter family of solutions. For $\det(M)=1$ (and $M \neq \pm I$)
the system has rank two and has  the one-parameter family of solutions
$(A, 2B, C) = \alpha(c, d-a , -b)$
for $\alpha \in \RR.$ 
Choosing $\alpha=1$, its discriminant is
$$\mbox{Disc}(Q):= (d-a)^2+4bc = Tr(M)^2 - 4.$$
Cases (1)-(3) correspond
to the discriminant being negative, positive or zero.

Diagonalizing $M$ over the complexe numbers
 gives either $\left[{x\atop 0}{0\atop
x^{-1}}\right]$ or $\pm\left[{1\atop 0}{x\atop 1}\right]$; in either
case, we find that any matrix $M'$ commuting with $M$ is of the form
$M'=eM+fI$.  Substituting this in, we find that $M'$ also preserves $Q$.
~~~$\bsq$ \\
%
%
\begin{theorem}\label{th55}
Suppose that $S_{ab}$ has irrational rotation number and that
 the $S_{ab}$ orbit of $(0,1)$
contains $(0,-1)$. Then the $T_{ab}$ orbit of 
$(0,1)$ contains $(0,-1)$, and 
$T_{ab}$ has a piecewise conic invariant circle.
If $S_{ab}^{(n)} (0,1) = (0,-1)$ for $n \in \ZZ$, 
then the number of conic pieces of the invariant circle is at most $|n|$,
and all conic pieces are of the same type, either arcs of ellipses,
arcs of hyperbolas or line segments, respectively.
\end{theorem}

\paragraph{Proof.}
The fact that the $T_{ab}$ orbit of 
$(0,1)$ contains $(0,-1)$ follows from Theorem~\ref{th34},
since irrational rotation number can only occur in case (ii).

Let $\sC$ be the closure of the orbit containing $\bv_0 :=(0,1)$;
we prove below that $\sC$ is an invariant circle with the 
required properties.
If so, it follows that all $T_{ab}$ orbits are scaled copies of this
invariant circle.

By hypothesis
 $S_{ab}^{(n)} (0,1) = (0,-1)$ with $n \in \ZZ$; the value of $n$
is unique since $S_{ab}$ has irrational rotation number.  Set $\sO=
\{S_{ab}^{(i)} (0,1) : 0 <i \le n \}$ if $n > 0$ and 
set $\sO = \{S_{ab}^{(i)} (0,1) : n <i \le 0 \}$ if $n < 0$;  
this is the minimal set $\sO$ to which 
Theorem \ref{th54} applies. It follows that the  rays determined by
the points of $\sO$ 
partition of the plane  into $|n|$ sectors, and 
on each sector the first return map of $T_{\mu\nu}$
is piecewise linear, with two
commuting pieces $M_1, M_2$. 

We consider the sector $J= \RR^{+}[\bv, \bv')$
in which $\bv$ is whichever of $(0,1)$ or $(0,-1)$ occurs
first in $\sO$. 
We have each $M_j \neq \pm I$ because, if not, then 
the rotation map $S_{ab}$ would
have a periodic point, and hence a rational rotation number,
a contradiction. 
Now Lemma~\ref{le51} applies, and shows that $M_1$ and $M_2$
each leave invariant a one-parameter family of quadratic forms,
multiples of $ Q_1$ and $Q_2$, respectively;
that is,  $M_1^T Q_1M_1 = Q_1$ and $M_2^T Q_2M_2 = Q_2.$
Moreover, since $M_1$ and $M_2$ commute, we can take $Q_1=Q_2$.

The claim  shows that
all the return visits of $\bv$ to $J$ lie on a single
level set of the quadratic form $Q_1$,
with the level set  parameter $\lambda$
specified by requiring that $\bv$ lie on this conic.
Since $S_{ab}$ has irrational rotation number the return visits
are dense on that part of this conic section in the sector $J$,
and the closure is the full  arc of the conic section, call it $\sC_1$,
that intersects the sector.
In the hyperbola case $|Tr(M_1)| > 2$, the sector $J$ cannot
include an eigenvector of $M_1$ (an asymptote), because if it
did, then the map $S_{ab}$ would 
have a periodic point, and hence a rational rotation number,
a contradiction; thus the sector $J$ contains an arc of 
one connected component of a hyperbola.
 
The arc $\sC_1$ is either an arc of an ellipse or a hyperbola, or
a line segment. 
By similar reasoning, on each other sector $J_i$, the visits of
the orbit of $\bv_0$ fill out the arc of a unique conic on
that sector. The image of the arc $\sC_1$, when iterated under
$T_{ab}$ will visit each sector; this visit is given by a
linear map, so the conic $\sC_i$  is the same type of conic (linear,
ellipse, or hyperbola) as $\sC_1$.

It remains to show that the union $\sC = \cup_{j=1}^n \sC_j$ 
of these pieces is a connected set. To see this, consider
the preimage $\bv_{-1}= T_{ab}^{-1}(\bv)$. This point is
not in $\sO$, so it necessarily
lies in the interior of some sector $J_i$, hence is contained in 
the conic piece $\sC_i$. Now, the forward
 iteration $T_{ab}(\sC_i)$
will cover the point $\bv_0$ in its interior. This image is
a conncected set, so the conic pieces $\sC$ and $\sC_j$ in the
sector $J_j$ adjacent to $J$, moving clockwise, must touch at
the point $\bv_0$. Next, forward iteration of $T_{ab}$ of a
connected neighborhood of $\bv$ over the full set $\sO$ 
shows that all segments
$\sC_j$ glue together into a connected set.
Thus $\sC$ is an invariant circle.~~~$\bsq$

\paragraph{Remark.}
The conclusion of
Theorem \ref{th55} need not hold without the requirement of irrational
rotation number, see Example~\ref{ex53} in \S4. \\
%
%
%
%
%
\section{Examples}
\hsp

In this section we  
first give  examples of one-parameter families where
$(0, 1)$ and $(0,-1)$ are in the same orbit, which
illustrate the theorems in \S3.
%
%

\begin{exam}\label{ex51}
If $1<a<\sqrt{2}$, and
$b=2{a^2-1\over a(a^2-2)}$,
then $T_{ab}$ is either periodic or has a piecewise elliptical
invariant circle. (These  parameter values come from assuming that
$T^{(8)}_{ab}(0,-1)=(0,1)$, with appropriate sign conditions.)
\end{exam}

A detailed analysis of  this example is given in the appendix,
where the following facts are proved.
One chooses $\sO= \{ T_{ab}^{(j)}: 1 \le j \le 8 \}$.
When the rotation number is rational, $T_{\mu\nu}$ is periodic;
when it is irrational Theorem~\ref{th55} applies to show that it 
has a piecewise conic invariant circle with at most eight pieces.
The conic sections are always ellipses. 

There is an explicit formula for the rotation number.
Set $a = 2\cos \theta$, with
$\pi/4 < \theta < \pi/3.$ 
The  rotation number  $r(S_{ab})$  is explicitly given by
$$ 
r(S_{ab}) = \frac{3\pi - 7\theta}{14 \pi - 32 \theta}.
$$
It follows that $r(S_{ab})$
 is rational if and only if $\theta$ is a rational
multiple of $\pi$.

\begin{figure}[htb]
\centerline{\psfig{file=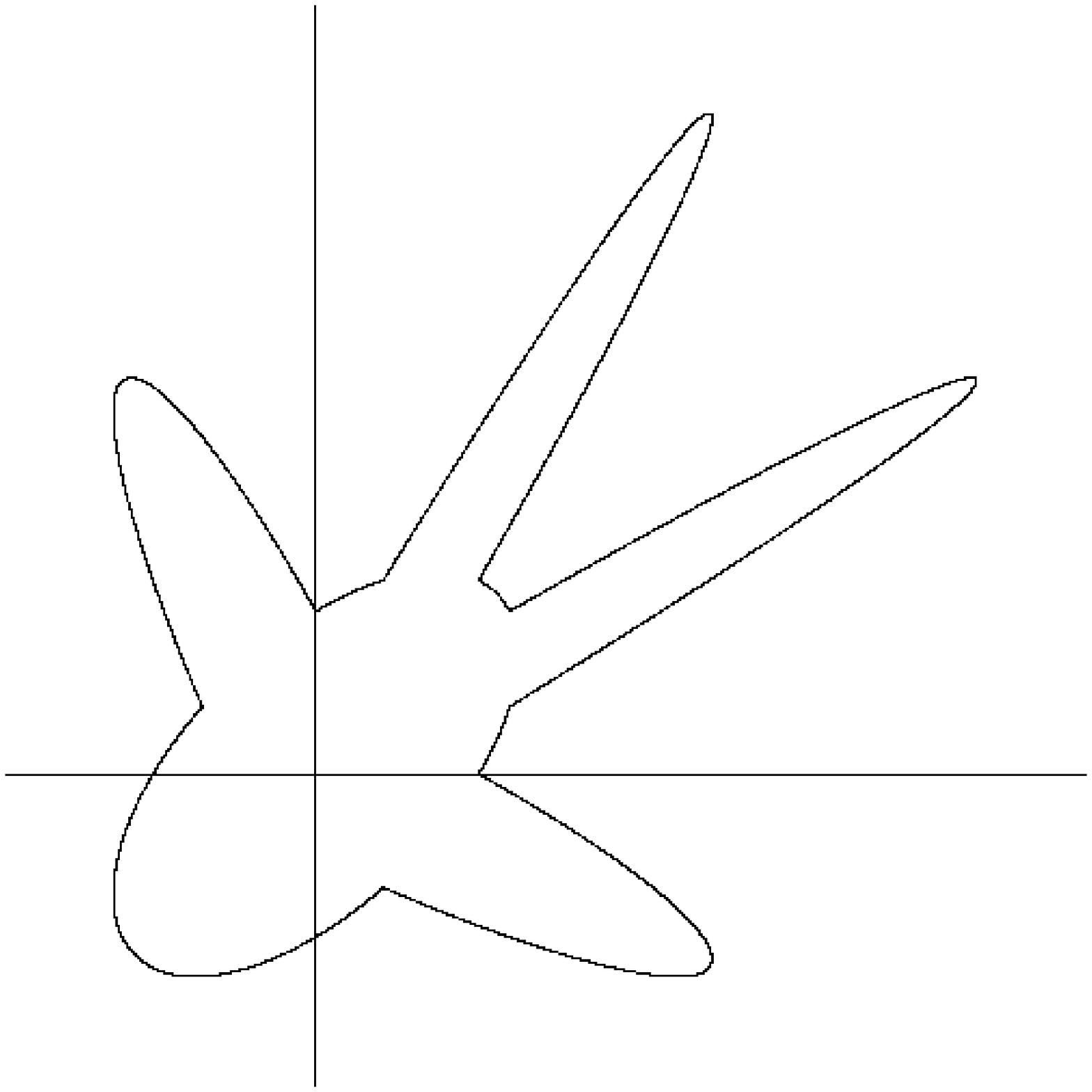,width=3in}}
\caption{Invariant circle for $a= {\root 4 \of 2}$, $b = - {\root 4 \of 2}$ 
 (``Bugs Bunny'').}
\label{fig1}
\end{figure}

The special case $\nu=0$, $\mu={\root 4 \of 2}$ 
(corresponding to  $a= -b = {\root 4 \of 2}$) 
has a rotation number which is provably irrational.
Indeed, one checks that $z=e^{i\theta}$ is a root of
$z^8 + 4z^6+4z^4 + 4z^2 + 1=0$; this equation
has no roots of unity as roots, so $\theta$ is not a rational
multiple of $\pi$.
An invariant circle of $T_{\mu\nu}$ is pictured in Figure \ref{fig1}; it is
a piecewise union of eight arcs of ellipses. (In this
case several of the arcs are parts of a single ellipse.)

A related example comes from  example 3.2 in part I,
which has  parameters
$a= \cos \frac{2\pi}{n}$ and $b = 2 \cos \theta$ with $0 < \theta < \pi$
for $\theta$ an irrational multiple of $\pi$, then
$r(S_{ab})= \frac{\theta}{\pi + n \theta}$ is irrational, and
this case has piecewise elliptic invariant circles.
%
%

%
%
\begin{exam}\label{ex52}
If $0<a<1$, and $b={1\over 2}{3a^2-3+\sqrt{a^4-2a^2+9}
\over a(a^2-2)}$, then $T_{ab}$ is either periodic or has a piecewise
conic invariant circle. All three cases of piecewise
ellipses, straight lines
or hyperbolas  occur, depending on the parameter value.
(These parameter values come from assuming that 
$T^{(10)}_{ab}(0,1)=(0,-1)$,
with appropriate sign conditions.)
\end{exam}

A detailed analysis of  this example is given in the appendix,
where the following facts are proved. We take
$\sO= \{ T_{ab}^{(j)}(0,1): 1 \le j \le 10 \}$.
In the case of rational
rotation number the map $T_{ab}$ is periodic.
For irrational rotation number  Theorem~\ref{th55}
applies, and there  are three cases, depending on the
trace of $M_1$. One has $0 <Tr(M_1) <2 $ for $1 < a < \alpha_0$
where 
$\alpha_0 \approx 0.3802775690976\ldots$ is
the real root of $x^4+3x^3+3x^2+x-1$ in the unit interval; 
invariant circles are then piecewise arcs of ellipses.
One has $Tr(M_1) = 2$ for $a = \alpha_0$; the invariant circle in
this case
is piecewise linear. One has $Tr(M_1) > 2$ for $\alpha_0 < a < 1$;
invariant circles are then  piecewise arcs of  hyperbolas.

\begin{figure}[htb]
\centerline{\psfig{file=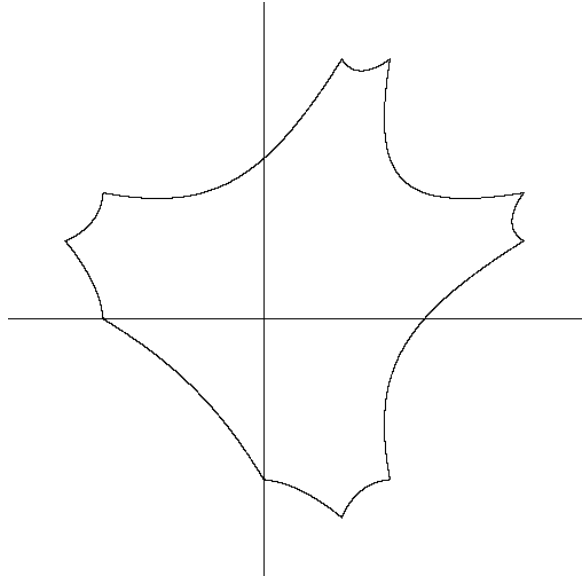,width=3in}}
\caption{Invariant circle for 
$a = \left( \frac{\sqrt{5} -1}{2} \right)^{1/2}$, 
$b= - \left( \frac{\sqrt{5}-1}{2} \right)^{1/2}$ (``Elmer Fudd'').}
\label{fig2}
\end{figure}

The following three figures picture
 invariant circles for  cases of Example \ref{ex52}.
Figure ~\ref{fig2} 
pictures the case $a=\sqrt{\sqrt{5}-1\over 2})\approx  0.78615 $,
which is the case  $b = -a$. The rotation number is
irrational (in fact, transcendental) and the invariant circle consists of 
ten arcs of hyperbolas. 

%
%

\begin{figure}[htb]
\centerline{\psfig{file=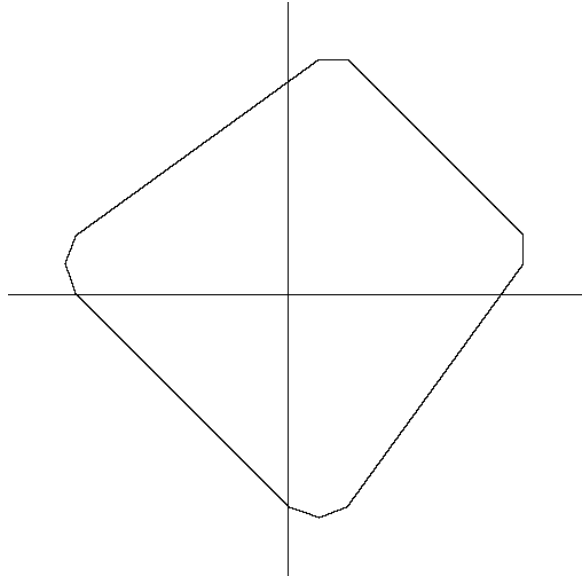,width=3in}}
\caption{Piecewise linear invariant circle for critical value
$a= \alpha_0 =0.380277$, 
$b$ as in example \protect\ref{ex52}.}
\label{fig3}
\end{figure}

Figure~\ref{fig3}  pictures the case of
the critical value
$a=\alpha_0 \approx 0.3802 $. The rotation number 
$$
r(T_{\mu\nu}) = \frac{2 \alpha_0^2 + 1}{9 \alpha_0^2 + 4}
$$
is an irrational algebraic number. Theorem~\ref{th55} applies
to show it has an invariant circle, which is piecewise linear
with ten pieces.

\begin{figure}[htb]
\centerline{\psfig{file=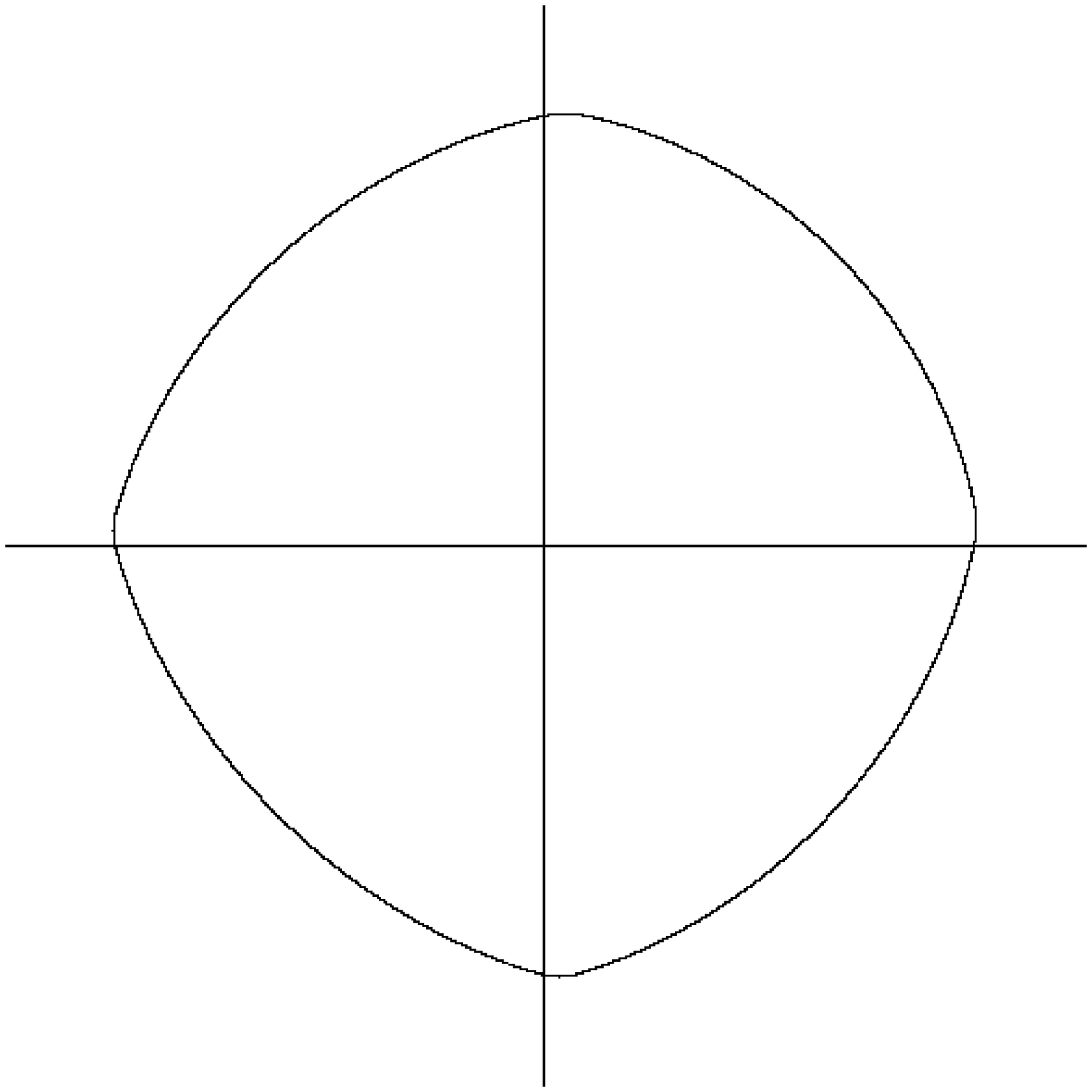,width=3in}}
\caption{Invariant circle for $a=0.1$, $b$ as in example \ref{ex52}.}
\label{fig4}
\end{figure}

Figure~\ref{fig4} pictures the case $\alpha=1/10$.
In this case the rotation number is irrational  
and the invariant circle consists of ten arcs of ellipses.

%
%

%
%

%
%
\begin{exam}\label{ex53}
If $1<a<\sqrt{2}$, and
$b=\frac{(a-1)(2a^3- 4a - 1)}{ a(a^2-2)(a^2 -a -1)}$,
then $T^{(13)}_{ab}(0,-1)=(0,1)$, but $S_{ab}$
has constant rational rotation number $r(S_{ab}) = \frac{1}{5}$
and $T_{ab}$
has no invariant circles. In this case all orbits diverge.
\end{exam}

The range of $b$ is $-\infty < b < 0$.
Theorem~\ref{th54} and \ref{th55} do not
apply, because all 
$S_{ab}$ in this family have constant rational 
rotation number $r=\frac{1}{5}$.
The following facts are proved in the appendix.
If we take $\sO = \{ T^{(j)}(0, -1) :~ 1 \le j \le 13 \}$
then part of the conclusion of Theorem~\ref{th54} still applies:
there exists a sector $J=\RR^{+}[\bv, \bv')$ given by
two elements of $\sO$  which has a
well defined first return map, which consists of exactly
two pieces, and the resulting matrices $M_1$ and $M_2$
commute. (However not all sectors have a well-defined first
return map.)
For the sector $J$ the  matrices $M_j$ have $|Tr(M_j)| > 2$, so their
invariant sets restricted to $J$ are hyperbolic arcs; 
however the asymptotes of the hyperbolas also lie in $J$.
These maps
fall under case (ii) of Theorem~\ref{th24}; all orbits
are unbounded.
Each $T_{ab}$ in this family satisfy
the following  weak form of the invariant circle
property:  
Each  orbit of $T_{ab}$ 
either falls on a finite number of rays, or else lies on
a finite union of (unbounded) hyperbolic arcs. 

A specific example is 
the parameter values $a = -b \approx 1.235877977$
where $a$ is a  root of $x^6-x^5-x^4-2x^2+3x+1=0.$ \\

It appears plausible that there are many more 
parameter values having invariant circles.
We present several non-rigorous numerical examples
showing numerical plots of orbits for other 
parameter values; these
are intended to  illustrate the structure
of (possible) general invariant circles.
The parameter values exhibited here were found by iterative
trial-and-error search, in which the first 10000 iterates 
appeared to form an invariant circle.
These include the values $(a,b)$ used in making the 
Figures \ref{fig5}--\ref{fig9} below.
One should view these figures 
with a grain of salt; 
a sufficiently slowly divergent periodic orbit could produce
a picture that looked like an invariant circle; qualitatively, the longer
the period of the periodic orbit, the closer the corresponding eigenvalue
seems to tend to be to 1, so an orbit with very high period will tend
to diverge very slowly, and will, furthermore, look like a curve.
However Conjecture A together with  various theorems in \S2 
would imply that
the set of parameter values having an irrational rotation
number has positive two-dimensional Lebesgue measure.
If so, this  makes it plausible that the computer plots actually
approximate invariant circles. That is,  
even if the parameter values below
don't correspond to invariant circles (as seems likely), 
then very small parameter changes should give some $T_{ab}$
having an invariant circle exhibiting the same qualitative appearance.
One may compare the apparant invariant circle in Figure \ref{fig9}
with one pictured in Froeschl\'{e} \cite{Fr68}, which stimulated
the work of Herman \cite[Chap. VIII]{He86}
on existence of invariant circles for these maps.
%
%
\begin{exam}
\label{ex54}
{\rm $a=0.2$, $b=-0.7$. (Figure \ref{fig5})  
At this scale the
invariant circle appears to be smooth; it seems unlikely that
this is actually the case.
}
\end{exam}
%
%
\begin{figure}[htb]
\centerline{\psfig{file=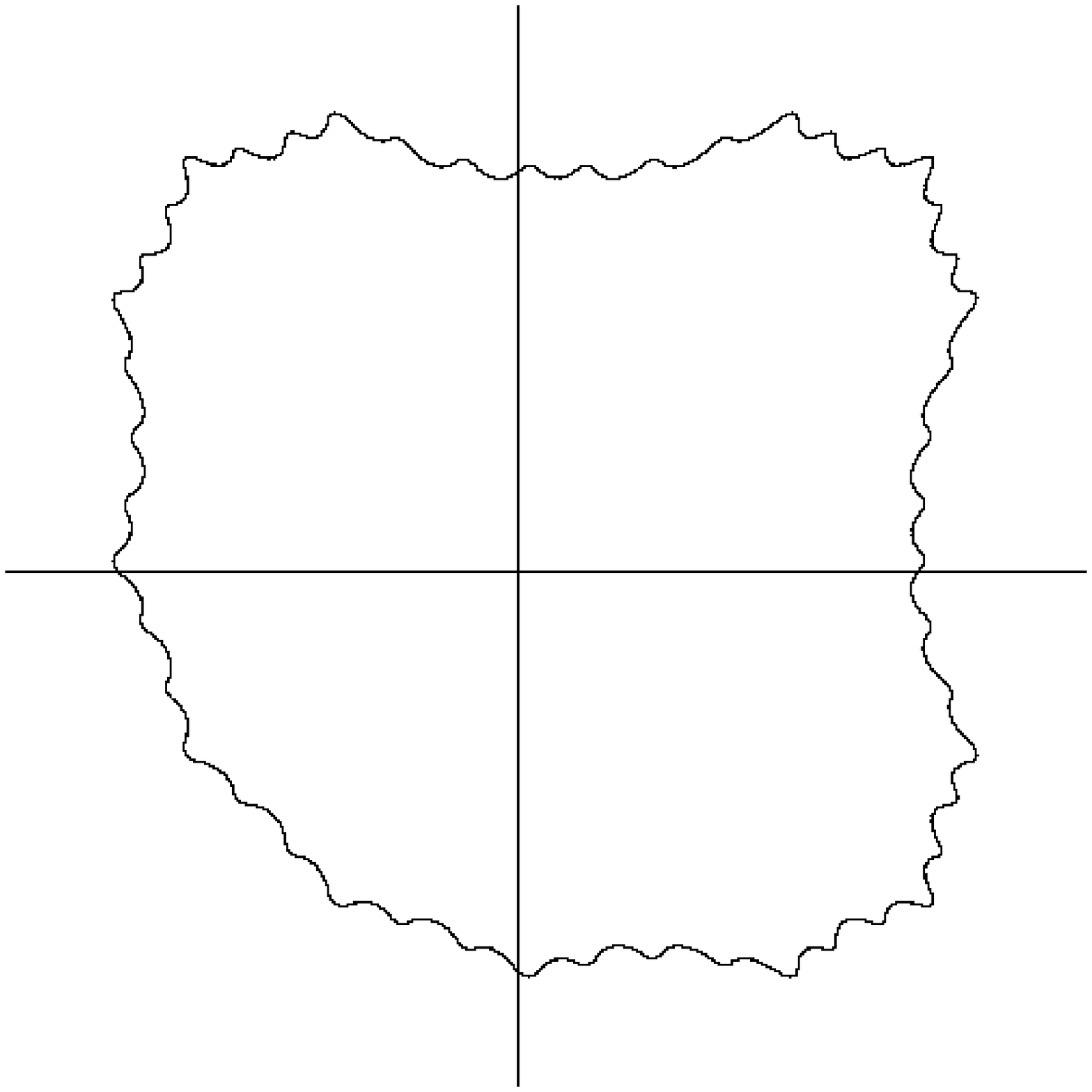,width=3in}}
\caption{Apparent invariant circle for $a=0.2$, $b=-0.7$.}
\label{fig5}
\end{figure}
%
%
\begin{exam}
\label{ex55}
{\rm
$a=1.4$, $b=-1.4$. (Figure \ref{fig6}) This
highlights a typical behavior of the generic invariant circle; there are
oscillations on several frequency scales (the thickness of the curve
results from oscillations of wavelength smaller than the pixel size).  The
``hair'' is usually not as pronounced as this.
}
\end{exam}
%
%
\begin{figure}[htb]
\centerline{\psfig{file=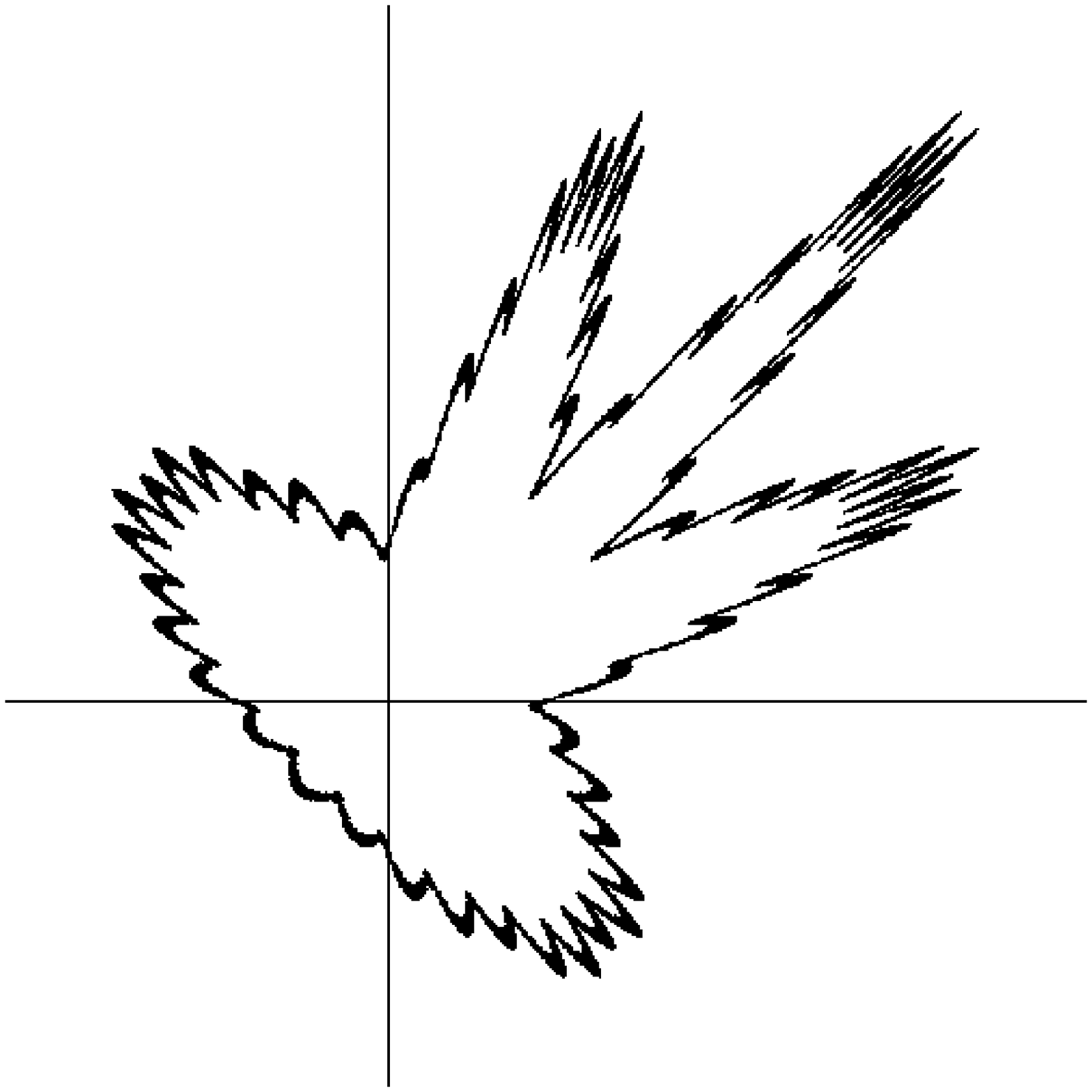,width=3in}}
\caption{Apparent invariant circle for $a=1.4$, $b=-1.4$.}
\label{fig6}
\end{figure}
%
%
\begin{exam}\label{ex56}
{\rm
$a=-.9$, $b=-4$. (Figure \ref{fig7}) The long
spikes in the second and fourth quadrants are typical of invariant circles
when $a,b<0$, although, again, not usually as pronounced, for larger
values of $b$.
}
\end{exam}

%
%
\begin{figure}[htb]
\centerline{\psfig{file=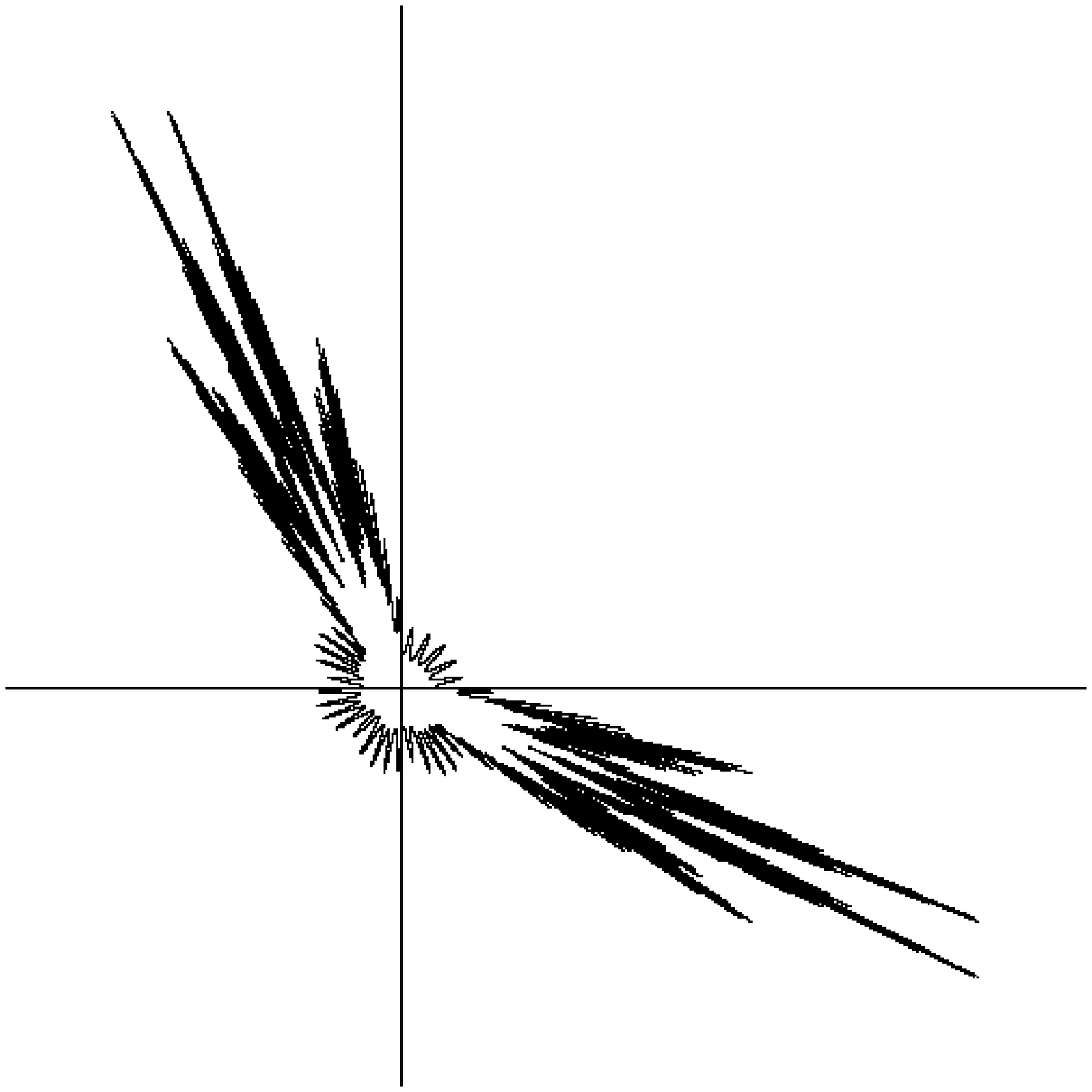,width=3in}}
\caption{Apparent invariant circle for $a=-0.9$, $b=-4$.}
\label{fig7}
\end{figure}

%
%
\begin{exam}\label{ex57}
{\rm
$a=1.5$, $b=1.1$. (Figure \ref{fig8}) For
$a,b>0$, on the other hand, the first and third quadrants features tend to
be pronounced.
}
\end{exam}
%
%
\begin{figure}[htb]
\centerline{\psfig{file=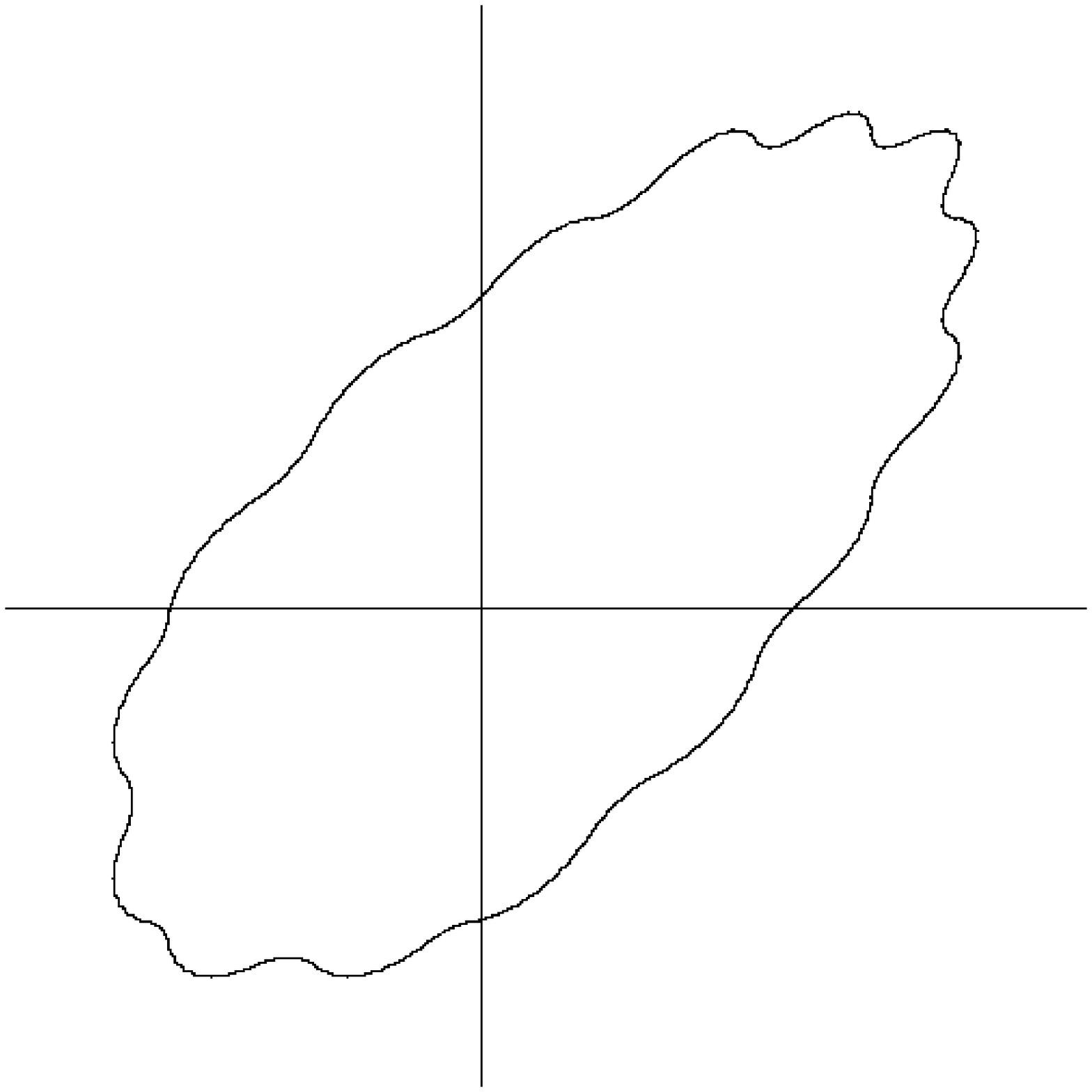,width=3in}}
\caption{Apparent invariant circle for $a=1.5$, $b=1.1$.}
\label{fig8}
\end{figure}
%
%

\begin{exam}\label{ex58}
{\rm
$a=1.9$, $b=-.2$. (Figure \ref{fig9}) For
$a>0>b$, the first
quadrant alone tends to be pronounced if $a>>-b$ (otherwise, 
as in figure \ref{fig5}),
no quadrant seems to dominate).
}
\end{exam}
%
%
\begin{figure}[htb]
\centerline{\psfig{file=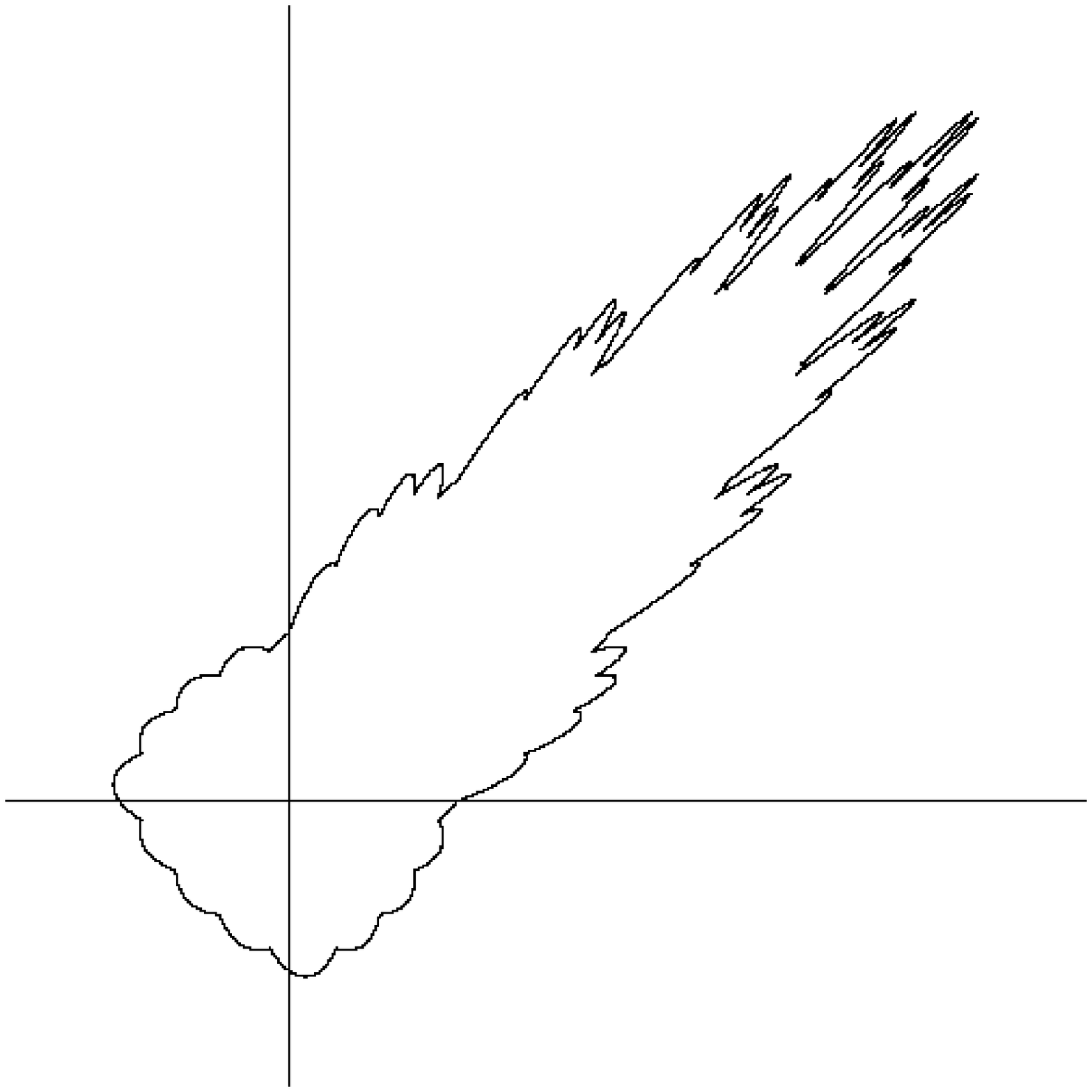,width=3in}}
\caption{Apparent invariant circle for $a=1.9$, $b=-0.2$.}
\label{fig9}
\end{figure}

%
%
%
%
%
\section{The Invariant Conic Set $\Omega_{Q}$ and $\Omega_{SB}$}

We  defined 
$\Omega_{Q}$ to be the set of parameter values $(a,b)$ 
for which $T_{ab}$ has a piecewise conic invariant circle of the type 
given in Theorem ~\ref{th54},
i.e. parameter values such that
the $S_{ab}$-orbit of $(0,1)$ contains $(0,-1)$, and $S_{ab}$ has irrational
rotation number. The fact that that the orbit of $(0,1)$ is
restricted makes this a rather small set.
%
%
\begin{theorem}\label{th72}
The set $\Omega_{P} \cup \Omega_{Q}$ has Hausdorff dimension $1$.
\end{theorem}

\paragraph{Proof.}
We have exhibited parametric curves in this set in the examples, so
it remains to show that $\Omega_{P}$ and $\Omega_{Q}$ each
have Hausdorff dimension at most one. 

We treat $\Omega_{Q}$, the argument for $\Omega_{P}$
being similar.
Theorem \ref{th54} gives the necessary condition
that either $S_{ab}^{(k)} (0,1) = (0,-1)$ or
$S_{ab}^{(k)} (0, -1) = S_{ab} (0,1)$, for some $k \in \ZZ_{> 0}$.
These conditions (namely $T_{ab}^{(k)}(0,1)_x = 0$, subject
to the constraint  $T_{ab}^{(k)}(0,1)_y < 0$, for each 
$k \in \ZZ \setminus \{0\}$)
cut out a countable number of pieces of real algebraic 
curves.  The resulting set
contains $\Omega_{Q}$.
It follows that $\Omega_{Q}$ has 
Hausdorff dimension at most $1$.~~~$\bsq$

In the introduction we proposed  the following conjecture.

\paragraph{Conjecture A.}
{\em The set $\Omega_{P} \cup \Omega_{Q}$ is a dense subset of 
$\Omega_{SB}$.} \\

We give  some justification for this conjecture.
For any positive $\epsilon$, if $r(S_{\mu\nu})$ is irrational, 
there is
a perturbation of $\nu$ smaller than $\epsilon$ such that $T_{\mu\nu'}$
maps $(0,1)$ to $(0,-1)$; that is, 
the perturbation goes to  $(0, - \lambda)$
for some positive $\lambda$ and  by Theorem 3.2 (ii) of part I,
it must be the case that  $\lambda=-1$.
One can rigorously show, for any rational, there is such a
point having rotation number between that rational and $r(S_{\mu\nu})$.
Experimentation suggests that this point will usually be in
$\Omega_{P}\cup\Omega_{Q}$, but the example \ref{ex53}
shows one way this can fail.

As an example, we consider the problem of finding points of
$\Omega_{per} \cup \Omega_{conic}$ close to the parameter values
$(a,b) = (1.4, -1.4)$ which has the (apparent) invariant circle 
pictured in Figure \ref{fig6}.
We found numerically that
$$
a=-b=1.399999999996869855138895062694261779494
$$
has a piecewise elliptic invariant circle (or is periodic), with
$T^{(810096)}(0,-1)=(0,1)$, complex eigenvalues. We also found numerically 
that
$$
a=-b=1.400000000000000366642660387254241957986
$$
has a piecewise hyperbolic invariant circle (or is periodic), with 
$T^{(35676178)}(0,1)=(0,-1)$, no eigenvectors in interval.

Both these parameter sets produce plots basically indistinguishable 
from Figure~\ref{fig6}. Since the simpler case has 810096 conic pieces, 
the piecewise conic structure is invisible in the computer plot. 

In fact, $\Omega_{P}\cup \Omega_{Q}$ is contained
in a countable union of
real algebraic curves; the argument above would say that any 2-dimensional
ball centered at a point of $\Omega_{SB}$
would contain a segment of such a curve.  In nearly all cases, 
as exhibited in Examples \ref{ex51} and \ref{ex52}, it appears that
such a curve does not have constant rotation number; if so, then it
necessarily contains a dense set of points with irrational rotation number,
and one can guarantee
a point of $\Omega_{conic}$ in any neighborhood, 
so that Conjecture A would follow. 

To conclude, we remark on the general question of the
existence of invariant circles when $r(S_{ab})$ has irrational
rotation number. 
A heuristic from KAM theory \cite{La93} is that
if the rotation number has very good rational approximations, then the
orbit has the possiblity to diverge, since it is very close 
for a long time to a number of
divergent parameter values.
On the other hand, if the rotation number
$r$ is irrational,
then the $S_{\mu\nu}$ orbit of $(0,1)$ will approach arbitrarily close to 
$(0,-1)$ infinitely often; this ought to place some constraints on the
extent to which the size of the iterates of $(0,1)$ can vary, 
and suggests the heuristic
that  under some Diophantine condition on $r$ this might
lead to an invariant circle. Indeed the result of 
Herman ~\cite[VIII.5.2]{He86}
shows this is the case when the rotation number is of constant
type (i.e. a badly approximable number). It remains an open
problem to decide if there exist parameters $(a,b)$ for
which $S_{ab}$ has irrational rotation number and $T_{ab}$
has no invariant circle. The  property of having no invariant
circle, if it occurs, 
will not come from a Diophantine property alone
of the rotation number $r(S_{ab})$; the particular parameter
values $(a, b)$ will matter. This is indicated by 
Example \ref{ex51}, which 
gives a family with continuously varying rotation number
in which all parameter values having irrational rotation number
necessarily have an invariant circle; all parameter values
in some open interval occur.
It might be helpful in resolving this open problem to
better understand the  mechanism by which
these examples, which can have 
an  irrational rotation number 
with arbitrary Diophantine properties, 
circumvent the KAM heuristic above.

%
%
%
%
%

\newpage
\section{Appendix: Rigorous Verification of Examples} 

In this appendix we give proofs for the assertions
made in \S5 about Examples~\ref{ex51}- \ref{ex53}.
The calculations for Example~\ref{ex51} below can be verified
by hand, but those for Example~\ref{ex52} and \ref{ex53}
require a computer.
The packages MAPLE and MAGMA were used for some of the
calculations below.
%
%
\noindent\paragraph{Example~\ref{ex51}.} {\em $T_{ab}$ with 
$1 < a < \sqrt{2}$, with
$b = 2 \frac{a^2 -1}{a^3 - 2a}.$ Here $b \in (-\infty, 0),$
and $T_{ab}^{(8)}(0, -1) = (0,1).$ }

\noindent\paragraph{Proof.} We first verify
that one has $T_{ab}^{(8)}(0, -1) = (0, 1).$ Set $\bv_0=(0, -1)$
and $\bv_j= T_{ab}^{(j)}(\bv_0)$; the sequence
of signs $S_j = \mbox{Sign}(\bv_j)_x)$ determine the
sequence of iterates. The resulting $\bv_j$ and their sign values
are given in Table~\ref{tableA1}.

\begin{table}[hp]\centering
\begin{tabular}{|l|r|}
\hline
\multicolumn{1}{|c|}{$\mbox{vector}$   } &
\multicolumn{1}{c|}{$\mbox{Sign}(\bv_j)_x)$  } 
\\ \hline
$\bv_1 = (1, 0)$              & $+$   \\
$\bv_2 = (a, 1)$              & $+$    \\
$\bv_3 = (a^2 -1, a)$         & $+$    \\
$\bv_4 = (a^3-2a, a^2 - 1)$   & $-$    \\
$\bv_5 = (a^2 -1, a^3 - 2a)$  & $+$    \\
$\bv_6 = (a, a^2 - 1)$        & $+$    \\
$\bv_7 = (1, a)$              & $+$    \\
$\bv_8 = (0, 1)$              & $\pm$  \\
\hline
\end{tabular}
\caption{Orbit $\sO$ for example \ref{ex51}}~\label{tableA1}
\end{table}

To verify the entries in Table~\ref{tableA1}, note that the condition
$1 < a < \sqrt{2}$ determines $\{\bv_j : 1 \le j \le 4\}$ and their
signs, and in particular
$$
\bv_5=( (a^3 - 2a)b + 1 - a^2, a^3 - 2a.) 
$$
We impose the condition that $\bv_5= R(\bv_4)$, and this requires
that
\beql{eqA01}
(a^3 - 2a)b + 2 - 2a^2 = 0,
\eeq
which specifies $b$ uniquely. Now the relation 
$T_{ab}^{-1}= R \circ T_{ab} \circ R^{-1}$ yields
$\bv_j = R(\bv_{9-j})$ for $5 \le j \le 8$, and also determines
their symbol sequence. This completes Table~\ref{tableA1}, and
shows $T_{ab}^{8}(0, -1) = (0,1)$.

We choose $\sO= \{\bv_j:~ 1 \le j \le 8\}$, and the hypotheses
of Theorem~\ref{th54} and \ref{th55} apply whenever the
rotation number $r(S_{ab})$ is irrational. The elements of $\sO$
appear in counterclockwise angular order as
$$ \bv_1, \bv_6, \bv_2, \bv_7, \bv_3, \bv_8, \bv_4,\bv_5. $$
We choose the sector $J := \RR^{+}[\bv, \bv') = \RR^{+}[\bv_4, \bv_5)$.
This sector  includes the entire third quadrant with signs $(-, -)$;
here $\bv_4$ and $\bv_5$ have signs $(-, +)$ and  $(+, -)$, respectively.
The breakpoint for $J$ is the ray determined by
\beql{eqA02}
\bv^{*}:= (0, -1).
\eeq
The first return maps on the subsectors $J_1= \RR^{+}[\bv, \bv^{*})$
and $J_2=\RR^{+}[\bv^{*}, \bv')$ are well-defined regardless of
whether the rotation number is rational or irrational.

The first return map on $J_1$ is $T_{ab}^{(5)}$ and has symbol
seqence $(-++++)$; its associated matrix\footnote
{To express the entries of $M_1$ in terms of the variable $a$ alone,
the variable $b$ is eliminated using \eqn{eqA01}.}
$M_1$ is:
\beql{eqA03}
M_1= \left[ {{a}\atop{a^2-1}}~~ {{-a^2+1}\atop{-a^3+2a}} \right].
\eeq
The first return map on $J_2$ is $T_{ab}^{(4)}$ 
with symbol sequence $(++++)$ and its associated
matrix $M_2$ is:
\beql{eqA04}
M_2 = \left[ {a^4 - 3a^2 + 1 \atop a^3 - 2a} ~~
{-a^3 + 2a \atop -a^2 + 1} \right].
\eeq
We have $Tr(M_1) = -a^3 +3a$ so
\beql{eqA05}
0 < Tr(M_1) < 2 \qquad \mbox{for}\quad 1 < a < \sqrt{2}.
\eeq

We parametrize the eigenvalues and eigenvectors of $M_1$ 
and $M_2$ using the parameter $\theta$, with
$$ 
a = 2 \cos \theta \quad\mbox{for}\quad 
\frac{\pi}{4} < \theta < \frac{\pi}{3}.
$$
The common eigenvectors $\bw^{+}$ and $\bw^{-}$ of $M_1$ and $M_2$  
are given by
$$ 
\bw^{\pm} = (e^{\pm i \theta}, 1). 
$$
A calculation (using $a = e^{i\theta} + e^{-i\theta}$) gives
$$
 M_1\bw_{+} = -e^{3i\theta}\bw^{+}, ~~~  M_1\bw^{-} = -e^{3i\theta}\bw^{-}.
$$
$$
M_2\bw^{+} = e^{4i \theta}\bw^{+}, ~~~  M_2\bw^{-}  = e^{-4i\theta}\bw^{-}.
$$
The matrices $M_1$ and $M_2$ leave invariant a common family
of ellipses centered at the origin. Such a family of  
ellipses centered at the
origin is specified by a single parameter $\phi$ with $0 \le \phi < 2\pi$,
and the family elements parametrized by $r>0$, with
$$
E= E(r, \phi) := 
\{ \bv(\omega)=( r \cos(\phi + \omega), r\cos(\omega)): 
0 \le \omega < 2\pi\}.
$$
The action of $M_1$ on
that portion of the ellipse in $J$ is to 
take $\bv(\omega)$ to $\bv(\omega + \pi - 3 \theta)$, while the
action of $M_2$ there 
takes $\bv(\omega)$ to $\bv(\omega -2\pi + 4\theta)$;
the multiples of $\pi$ are specified to keep the
image angle in $[0, 2 \pi]$.

We now show that the rotation number $r(S_{ab})$ is given by
\beql{eqA08}
r(S_{ab})= \frac{3 \pi - 7 \theta}{14\pi - 32 \theta}.
\eeq
To establish this, we consider iteration $T_J^{(p)}$ of the
first return map $T_J$ on the sector $J$. Pick $\bw \in J$
lying on $E(r, \phi)$ 
and set
$$ T_J^{(p)}(\bw) = M_1^{m_1}M_2^{m_2}\bw ,$$
with $m_1= m_1(p), m_2= m_2(p),$ both going to $+\infty$
as $p \to \infty$. Since each appearance of $M_1$ or $M_2$ 
corresponds to one counterclockwise revolution of the circle,
while $M_1$ and $M_2$ correspond to $5$ and $4$ iterations
of $T_{ab}$, respectively, we have
\beql{eqA09}
r(S_{ab}) = \lim_{p \to \infty} \frac{m_1(p) + m_2(p)}{5m_1(p) + 4 m_2(p)}.
\eeq
However the  condition that the angle $\omega$ remains inside
the sector $J$ of $E(r, \phi)$ at each iteration of $T_J$
yields
$$
0 < \phi + m_1(p)(\pi - 3\theta) + m_2(p)(-2\pi + 4 \theta) < 2\pi.
$$
Letting $p \to \infty$, this constraint implies that
\beql{eqA10}
\lim_{p \to \infty} \frac{m_2(p)}{m_1(p)}= 
\frac{\pi - 3 \theta}{2\pi - 4\theta}.
\eeq
Rewriting  \eqn{eqA09} as
$$
r(S_{ab}) = \lim_{p \to \infty}~  \frac{ 1+ \frac{m_2(p)}{m_1(p)}}
{5 + 4 \frac{m_2(p)}{m_1(p)}}
$$
and substituting \eqn{eqA10} yields the desired 
rotation number formula \eqn{eqA08}.

Suppose that $r(S_{ab})$ is rational.  By \eqn{eqA08} this
occurs if and only if $\theta$ is a rational multiple of
$\pi$.
We show that $T_{ab}$ is periodic.
Indeed, for rational rotation number $S_{ab}$ has a periodic
orbit, and this orbit necessarily contains some point $\bw \in J$.
Since $S_{ab}^{(q)}(\bw) = \bw$ where $q$ is the period, we have
$T_{ab}^{(q)}(\bw) = \lambda \bw$ for some real $\lambda$. 
Since this value is an iterate of a first return map to $J$, there
is some $p$ such that 
$$T_{ab}^{(q)}(\bw) = T_J^{(p)}(\bw)= M_1^{m_1}M_2^{m_2}\bw . $$ Now
$M_1^{m_1}M_2^{m_2}$ has $\bw$ as an eigenvector with a real eigenvalue,
but it also has $\bw^{+}$ and $\bw^{-}$ as eigenvectors with complex
conjugate eigenvalues of norm 1. It follows that $M_1^{m_1}M_2^{m_2}= \pm
I$.  so that
$$T_{ab}^{(2q)}(\bw) = T_J^{(2r)}(\bw) = M_1^{2m_1}M_2^{2m_2}\bw.$$
It follows that  the first return map $T_J^{(2p)}$ is the 
identity on a subsector of $J$ of positive width, so $S_{ab}$ must 
have infinitely many periodic points of period $2q$, and
Theorem~\ref{th24}(iii) shows that $T_{ab}$ is periodic,
with  period dividing $2q$.

Suppose that $r(S_{ab})$ is irrational. Then Theorem~\ref{th55}
applies to show that $T_{ab}$ has a piecewise conic invariant
circle with at most $8$ pieces, and the pieces are arcs of
ellipses by \eqn{eqA05}.

 We consider the special parameters $a= -b = {\root 4 \of 2}$.
Defining $\theta$ by  $a = 2 \cos \theta$ shows that
$x= e^{i \theta}$ satisfies $x + 1/x ={\root 4 \of 2},$
hence it also satisfies
\beql{eqA11}
 x^8 + 4 x^6 + 4x^4 + 4x^2 + 1 = 0. 
\eeq 
We show the rotation number $r(S_{ab})$ is
irrational by showing that $\frac{\theta}{\pi}$ is irrational,
i.e. that $e^{i\theta}$ is not a root of unity.
It suffices to show that \eqn{eqA11} has no root
that is a root of unity. One checks using MAPLE
that this polynomial is irreducible over $\QQ$, and that it is
not a cyclotomic polynomial since it has a root off the
unit circle. Thus it has no roots of unity as roots,
and $r(S_{ab})$ is irrational. Thus  $T_{ab}$ has a 
piecewise elliptical invariant circle, pictured in 
Figure~\ref{fig1}.~~~$\bsq$

%
%

\noindent\paragraph{Example~\ref{ex52}.} {\em $T_{ab}$ with 
$0 < a < 1$, with
$b = 2 \frac{a^2 -1}{a^3 - 2a}.$ Here $T_{ab}^{(10)}(0, 1) = (0,-1).$ }

\noindent\paragraph{Proof.} We first verify
that one has $T_{ab}^{(10)}(0, 1) = (0, -1).$
Set $\bv_0=(0, -1)$
and $\bv_j= T_{ab}^{(j)}(\bv_0) = \bv_j,$; the 
signs  $S_j = \mbox{Sign}(\bv_j)_x)$ determine the
sequence of iterates. The resulting $\bv_j$ and their sign values
are given in Table~\ref{tableA2}.

\begin{table}[hp]\centering
\begin{tabular}{|l|r|}
\hline
\multicolumn{1}{|c|}{$\mbox{vector}$} &
\multicolumn{1}{c|}{$\mbox{Sign}(\bv_j)_x$} 
\\ \hline
$\bv_1 = (-1, 0)$                            & $-$    \\
$\bv_2 = (-b, -1)$                           & $+$    \\
$\bv_3 = (-ab+1, -b)$                        & $+$    \\
$\bv_4 = (-a^2b+ a+b , 1-ab)$                & $+$    \\
$\bv_5 = (a^3b +a^2 + 2ab -1, -a^2b + a+b)$  & $-$    \\
$\bv_6 = (-a^2b+a +b, a^3b +a^2+2ab -1)$     & $+$    \\
$\bv_7 = (-ab+1, -a^2b + a + b )$            & $+$    \\
$\bv_8 = (-b, -ab + 1)$                      & $+$    \\
$\bv_9 = (-1, -b) $                          & $-$    \\
$\bv_{10} = (0, -1)$                         & $\pm$  \\
\hline
\end{tabular}
\caption{Orbit $\sO$ for example \ref{ex52}}~\label{tableA2}
\end{table}

To verify the entries in Table~\ref{tableA2}, note that the condition
$0 < a < 1$ and the given sign sequence
determine $\{\bv_j : 1 \le j \le 5\}$, and in particular
$$
\bv_6=(a^3b^2 + 2a^2b + 2ab^2-a-2b,a^2 +2ab+a^3b -1) 
$$
We impose the condition that $\bv_6= R(\bv_5)$, and this requires
that
\beql{eqA12}
(a^3-2a)b^2 - 3(a^2-1)b + 2a=0.
\eeq
The correct choice of root of this quadratic equation determines
$b$ as given; one has $-\sqrt{2}< b < 0$ and 
 it produces the sign sequence
as given above for $1 \le j \le 5$. As in
Example~\ref{ex51}  the relation $\bv_6= R(\bv_5)$ forces 
$\bv_{j} = \bv_{11-j}$ for $6 \le j \le 10$, and determines
the remaining sign sequence, and gives
$T_{ab}^{(10)}(0, -1)= (0,1).$

We take $\sO= \{ \bv_j:~ 1 \le j \le 10 \}$,  and the hypotheses
of Theorem~\ref{th54} and \ref{th55} apply whenever the
rotation number $r(S_{ab})$ is irrational. The elements of $\sO$
appear in counterclockwise angular order as
$$ \bv_1, \bv_{10}, \bv_6, \bv_2, \bv_7, \bv_3,\bv_8, \bv_9, \bv_5. $$
We choose the sector $J=\RR^{+}[\bv, \bv') := \RR^{+}[\bv_1, \bv_{10})$,
which is exactly the third quadrant, with sign $(-,-).$ The breakpoint is
\beql{eqA13}
\bv^{*}= (a^2 -1, (a^2-1)b -a).
\eeq
The first return maps on the subsectors $J_1= \RR^{+}[\bv, \bv^{*})$
and $J_2=\RR^{+}[\bv^{*}, \bv')$ are well-defined regardless of
whether the rotation number is rational or irrational.

The first return map on $J_1$ is $T_{ab}^{(4)}$ 
with symbol sequence $(-++-)$ and its associated
matrix $M_1$ is:
$$
M_1 = \left[ {(a^2-1)b^2 -2ab +1 \atop (a^2-1)b - a}~~ 
{(1-a^2)b + a \atop 1-a^2} \right].
$$
The first return map on $J_2$ is $T_{ab}^{(13)}$ 
with symbol sequence $(-+++-+++-+++-)$ and its associated
matrix $M_2$ is:
$$
M_2 = \left[ {(1-a^2)b + a \atop 1-a^2}~~ 
{a^2-1 \atop \frac{a^2(a^2-2)}{(a^2-1)b-a}} \right].
$$
We have $Tr(M_1) = (a^2-1)(b^2-1) - 2ab + 1$ and using the
relation \eqn{eqA12} and the fact that $b<0$, we obtain
\footnote{Take the resultant of \eqn{eqA12} and the polynomial
$Tr(M_1)-2.$ It is a polynomial of degree $8$ in $a$ which has
two degree $4$ factors over $\QQ$. The other degree $4$ factor
produces extraneous roots.}
$$ 
Tr(M_1)=2 \quad\mbox{for}\quad a= \alpha_0 \approx 0.3802,
$$
where $\alpha_0$ is the unique root in the unit interval of
\beql{eqA15}
x^4 + 3x^3 + 3x^2 +x - 1=0.
\eeq
One can check that $0 < Tr(M_1)< 2$ for
$0 < a < \alpha_0$ and $Tr(M_1) > 2$ for $ \alpha_0 < a < 1$. \\

{\bf Subcase 1.} {\em Piecewise hyperbolic case ($\alpha_0 < a < 1$)} \\

The matrix $M_1$ is hyperbolic and has real eigenvectors
and eigenvalues; therefore $M_2$ does also, since it commutes
with $M_1$ and has the same eigenvectors. One checks
$Tr(M_2) > 2$ as well.
The first return map $T_J$ for $J$ leaves invariant hyperbolic
arcs in $J$, even in the rational rotation number case. These
are bounded arcs, for one can verify numerically that the 
asymptotes of the hyperbola fall strictly outside $J$. 
Thus the orbit of a given point $\bw \in J$,
lies on some  bounded hyperbolic arc $\sC_1$ inside $J$.

Let $\lambda_1= \lambda_1(a) >1$ be the large eigenvalue of $M_1$,
the other eigenvalue being $1/\lambda_1$; let $\lambda_2 >1$
be the large eigenvalue of $M_2$, the other being $1/\lambda_2$.
We choose the eigenvectors so that 
$M_1\bw_+= \lambda_1 \bw_+$
and $M_1\bw_{-}= 1/\lambda_1 \bw_{-}.$
We now claim that
\beql{eqA16}
 M_2\bw_{+}    = \frac{1}{\lambda_2}\bw_{+}.
\eeq
To establish this, choose a point $\bw$ in $J$, whose first
return iterates $T_J^{(p)}$ necessarily all
lie on a fixed hyperbolic arc $\sC_1$ in $J$. Write
$$ T_J^{(p)}(\bw) = M_1^{m_1}M_2^{m_2}\bw ,$$
with $m_1= m_1(p), m_2= m_2(p),$ both going to $+\infty$
as $p \to \infty$. Express $\bw$ in terms of the eigenvectors as
$$
\bw = c_1 \bw_{+} + c_2\bw_{-},
$$
with $c_1c_2 \ne 0$ since $\bw$ cannot be an eigenvector.
We now argue by contradiction, and suppose
$M_2 \bw_{+}= {\lambda_2}\bw_{+}.$ Then 
$$
T_J^{(p)}(\bw) = c_1 \lambda_1^{m_1}\lambda_2^{m_2} \bw_{+} +
c_2 \lambda_1^{-m_1}\lambda_2^{-m_2}\bw_{-},
$$
The first term on the right  gets large as $p \to \infty$, while
the second gets small, hence
 $||T_J^{(p)}(\bw)|| \to \infty$ as $p \to \infty$.
This contradicts  $T_J^{(p)}(\bw) \in \sC_1$ for all $p$; the claim follows.

We show that the rotation number $r(S_{ab})$ is given by
\beql{eqA17}
r(S_{ab}) = \frac{\log \lambda_2 + 3 \log \lambda_1}
{4 \log \lambda_2 + 13 \log \lambda_1}.
\eeq
We pick a point $\bw \in J$ and study the first return map iterates,
with $T_J^{(p)}(\bw) = M_1^{m_1}M_2^{m_2}\bw $ and
$m_1=m_1(p), m_2=m_2(p) \to \infty$ as $p \to \infty$. 
The return map for $J_1$ goes around the circle once, while that
for $M_2$ goes around three times, and using the number of
steps of $T_{ab}$ each takes we obtain
\beql{eqA18}
r(S_{ab}) = \lim_{p \to \infty} \frac{m_1(p) + 3m_2(p)}
{4m_1(p) + 13m_2(p)}.
\eeq
Expressing $\bw= c_1\bw_{+} + c_2\bw_{-}$
in terms of eigenvectors and using \eqn{eqA16}
yields 
$$
T_J^{(p)}(\bw) = c_1 \lambda_1^{m_1}\lambda_2^{-m_2}\bw_{+} +
  c_2 \lambda_1^{-m_1}\lambda_2^{m_2}\bw_{-}.
$$
In order for these vectors to remain on $\sC_1$ their norms
must be bounded away from $0$ and $\infty$, so there is a positive constant
$C$ such that $1/C< \lambda_1^{m_1}\lambda_2^{-m_2} < C.$
This yields
$$
\lim_{p \to \infty}~ \frac{m_1(p)}{m_2(p)} = 
\frac{\log \lambda_2}{\log \lambda_1}.
$$
Substituting this in \eqn{eqA18} yields the rotation number formula
\eqn{eqA17}.

When $r(S_{ab})$ is rational, then $T_{ab}$ is periodic.
To show this, we use the fact that the eigenvectors of $M_1$ and
$M_2$ fall on rays lying outside the sector $J$. If $r(S_{ab})$ is rational,
then $S_{ab}$ has a periodic orbit, which necessarily visits the sector
$J$. Choosing an element $\bw \in J$ on this orbit, one has
some first return value $T_J^{(p)}(\bw) = \lambda \bw,$
for some real value $\lambda$. Now 
$T_J^{(p)}(\bw) = M_1^{m_1}M_2^{m_2}\bw $,
and the matrix $M = M_1^{m_1}M_2^{m_2}$ now has three distinct 
eigenvectors, $\bw_{+},~ \bw_{-}$ and $\bw$, so $M= \pm I$.
Now we find that $T_J^{(2p)}(\bw') =  M_1^{2m_1}M_2^{2m_2}\bw'= \bw'$
for all $\bw'$ on a subsector of positive angular width.
It follows that $S_{ab}$ has infinitely many periodic orbits, so
$T_{ab}$ is periodic by Theorem~\ref{th24}(iii).

When $r(S_{ab})$ is irrational, Theorem~\ref{th55} applies to show
that $T_{ab}$ has an invariant circle which is a piecewise union
of at most $10$ arcs of hyperbolas.

We consider the special parameter values 
$a= -b= \sqrt{\frac{\sqrt{5}-1}{2}}.$
Here $Tr(M_1)= a^4 + 2 > 2.$ The larger eigenvalue $\lambda_1$ of $M_1$
satisfies $x^2 - (a^4+2) x + 1=0$,
and hence one finds it is an algebraic unit that satisfies
$$
x^4 - 7x^3 + 13x^2 -7x +1 =0.
$$
A similar calculation shows that $\lambda_2$ is an algebraic unit
that satisfies
$$
x^8 + 23x^6 - 77x^4 + 23x^2 + 1=0.
$$
The rotation number formula \eqn{eqA17} gives
\beql{eqA21}
r(S_{ab}) = \frac{\log (\lambda_1^3 \lambda_2)}
{\log(\lambda_1^{13}\lambda_2^4)},
\eeq 
and the condition for $r(S_{ab})$ to be irrational is that
$$ \lambda_1^{n_1}\lambda_2^{n_2} \neq 1 ~~\mbox{for}~~
(n_1, n_2) \in \ZZ^2 \backslash (0,0).
$$
That is, $\lambda_1$ and $\lambda_2$ must be multiplicatively
independent algebraic units in the algbraic number field over $\QQ$ that they
generate. This can be verified by checking that their logarithms and
the logarithms of the absolute values of their algebraic conjugates
span a two dimensional lattice. Given that $r(S_{ab})$ is irrational,
Theorem~\ref{th55} applies to show $T_{ab}$ has an invariant circle,
which consists of at most 10 arcs of hyperbolas, as pictured in
Figure~\ref{fig2}.

Given that $r(S_{ab})$ is irrational, we can show it is transcendental.
Clearing denominators in \eqn{eqA21}
and exponentiating gives
$$
(\lambda_1^{13}\lambda_2^{4})^r = \lambda_1^3 \lambda_2.
$$
The right side of this equation is an algebraic number, and
the left side is an algebraic number raised to the $r$-th power.
If $r= r(S_{ab})$ were an irrational algebraic number, then
 the left side
would be transcendental by the criterion
\footnote{In \cite[Theorem 2.4]{Ba75} take
 $\alpha_1= \lambda_1^{13} \lambda_2^{4}$ and $\beta= r$,
and $[1, \beta]$ are linearly independent over $\QQ$ since $r$
is irrational.}
of Theorem 2.4 in Baker~\cite{Ba75}. This is a contradiction, so
we conclude that $r(S_{ab})$ is transcendental. \\

{\bf Subcase 2.} {\em Piecewise linear  case ($a= \alpha_0$)} \\

Here $a = \alpha_0$, $b= \alpha_0^3 + 2\alpha_0^2 + \alpha_0 - 1$,
where $\alpha_0$ satisfies \eqn{eqA15}. 
The sector $J = \RR^{+}[(-1,0)), (0, -1))$ and in
 this case one finds that
$$
M_1 = I + \left[ {\alpha_0^2 \atop -\alpha_0^2}~~
{\alpha_0^2 \atop -\alpha_0^2} \right],
$$

$$
M_2 = I + \left[ {\alpha_0^2-1 \atop 1-\alpha_0^2}~~
{\alpha_0^2-1 \atop 1-\alpha_0^2} \right].
$$
Consequently
\beql{eqA24} 
M_1^{m_1}M_2^{m_2} = I +  \left[ {x \atop -x}~~{x \atop -x} \right],
\eeq
with $x = (m_1 + m_2)\alpha_0^2 - m_2.$ By Lemma~\ref{le51} these
commuting maps leave invariant a family of parallel straight lines,
and in the sector $J$ the first return map
iterates of any given point $\bw$ lie on
a certain finite line segment $\sC_1$.

We show that the rotation number
\beql{eqA25}
r(S_{ab})= \frac{2 \alpha_0^2 +1} {9\alpha_0^2 + 4}.
\eeq
For this, picking a point $\bw \in J$ and using the
fact that the orbit $T_J^{(p)}(\bw) = M_1^{m_1}M_2^{m_2}\bw $
Just as in subcase 1) we have 
\beql{eqA26}
r(S_{ab})= \lim_{p \to \infty}~ \frac{m_1(p)+3m_2(p)}{4m_1(p)+ 13m_2(p)}.
\eeq
Next, using the fact that the
orbit $T_J^{(p)}(\bw)$ lies on the line segment $\sC_1$, so
has norm bounded away from $0$ and $\infty$, using \eqn{eqA24} we have
$|m_1(p) \alpha_0^2 + m_2(p)(\alpha_0^2 - 1)| < C$ for some positive $C$.
Letting $p \to \infty$ yields
$$
\lim_{p \to \infty} \frac{m_1(p)}{m_2(p)}= \frac{1- \alpha_0^2}{\alpha_0^2},
$$
and substituting this in \eqn{eqA26} yields \eqn{eqA25}.

The rotation number $r(S_{ab})$ is irrational, as is verified by
computing the irreducible polynomial it satisfies over $\QQ$,
which has degree $4$. Now Theorem~\ref{th55} applies to show
that $T_{ab}$ has a piecewise linear invariant circle with at
most $10$ pieces. This justifies Figure~\ref{fig3}. \\

\noindent\paragraph{\bf Subcase 3.}
 {\em Piecewise ellipse case ($0 < a< \alpha_0$)} \\

This case is analyzed similarly to the case of Example~\ref{ex51}
and we omit the details. The invariant circle is piecewise conic
with arcs of ellipses when the rotation number is irrational,
and $T_{ab}$ is periodic when the rotation number is rational.

Consider the special parameter values $a = \frac{1}{10}$,
where $b = 5( \frac{297 - \sqrt{89801}}{199}).$ One finds that the
(complex conjugate) eigenvalues $\lambda_1, \bar{\lambda}_1$
of $M_1$ and $\lambda_2, \bar{\lambda}_2$ of $M_2$ are roots of
degree $4$ polynomials over $\ZZ[x]$ with large coefficients.
Verifying irrationality of the rotation number $r(S_{ab})$ 
reduces to showing that $<\lambda_1, \lambda_2>$ generate 
multiplicatively a
free abelian group of rank $2$. These numbers
are not algebraic integers, and it suffices to give a set of prime ideals
such that the valuations of $\lambda_1$ and $\lambda_2$ at those primes
span a lattice of  rank $2$.  
In particular, at the four prime ideals  above 2, $\lambda_1$
has valuations $(2,-2,2,-2)$, while $\lambda_2$ has valuations
$(7,-3,3,-7)$.  Thus $\lambda_1^i\lambda_2^j$ is an algebraic integer
only when $i=j=0$.
The irrationality of  $r(S_{ab})$
shows that this $T_{ab}$ has a piecewise elliptic invariant circle, as
pictured in Figure~\ref{fig4}. 
~~~$\bsq$

%
%
\noindent\paragraph{Example~\ref{ex53}.} 
{\em $T_{ab}$ with 
$1 < a < \sqrt{2}$, with
$b=\frac{(a-1)(2a^3- 4a - 1)}{ a(a^2-2)(a^2 - a - 1)}$. 
Here $T_{ab}^{(13)}(0, -1) = (0,1).$ }

\noindent\paragraph{Proof.}
One verifies 
that one has $T_{ab}^{(13)}(0, -1) = (0, 1)$ in a manner similar
to the earlier examples.
Set $\bv_0=(0, -1)$
and $\bv_j= T_{ab}^{(j)}(\bv_0),$ and define the sequence
of symbols $S_j = \mbox{Sign}(\bv_j)_x$ which determine the
sequence of iterates. The resulting $\bv_j$ and their symbol values
$S_j$ are given in Table~\ref{tableA3}. The value 
$b= \frac{(a-1)(2a^3- 4a - 1)}{ a(a^2-2)(a^2 - a - 1)}$  is determined by
the requirement that $R(\bv_7)= \bv_7$, plus the 
given sign seqence $\{S_i: 1 \le i \le 6\}.$ This value of $b$ has been
substituted in the earlier iterates; one verifies retrospectively
that the symbol sequence $S_j:~ 1 \le \bv_j \le 6\}$ is as indicated.


\begin{table}[hp]\centering
\begin{tabular}{|l|r|}
\hline
\multicolumn{1}{|c|}{$\mbox{vector}$} &
\multicolumn{1}{c|}{$\mbox{Sign}(\bv_j)_x$} 
\\ \hline
$\bv_1 = (1, 0)$                                          & $+$    \\
$\bv_2 = (a, 1)$                                          & $+$   \\
$\bv_3=  (a^2 -1, a)$                                      & $+$    \\
$\bv_4 = (a(a^2-2), a^2 - 1)$                             & $-$    \\
$\bv_5 = (a(a-1)(a^2-2)/(a^2-a-1), a(a^2 - 2)$            & $+$    \\
$\bv_6 = (a(a^2-2)/(a^2-a-1),a(a-1)(a^2-2)/(a^2-a-1))$    & $+$    \\
$\bv_7 = ( a(a^2-2)/(a^2-a-1),  a(a^2-2)/(a^2-a-1))$      & $+$    \\
$\bv_8 = (a(a-1)(a^2-2)/(a^2-a-1), a(a^2-2)/(a^2-a-1))$   & $+$    \\
$\bv_9 = ( a(a^2-2), a(a-1)(a^2-2)/(a^2-a-1)) $           & $-$    \\
$\bv_{10} = (a^2 - 1), a(a^2-2))$                         & $+$    \\
$\bv_{11} = (a, a^2 - 1)$                                 & $+$    \\
$\bv_{12} = (1, a))$                                      & $+$    \\
$\bv_{13} = (0, 1)$                                       & $\pm$  \\
\hline
\end{tabular}
\caption{Orbit $\sO$ for example \ref{ex53}}~\label{tableA3}
\end{table}

The points appear around the circle in counterclockwise order
$$ 
\bv_1, \bv_6, \bv_{11}, \bv_2, \bv_7, \bv_{12}, \bv_3, \bv_8, \bv_{13},
\bv_4, \bv_9, \bv_5, \bv_{10}.
$$
Of the thirteen arcs these determine, five have a well-defined
first return map and eight do not. The ones that do are
divided into two intervals by a breakpoint and the conclusion of
Theorem~\ref{th54} holds for them; they are the
intervals with left endpoint 
$\bv_9, \bv_{10}, \bv_{11}, \bv_{12}, \bv_{13}$, respectively.
The other
eight intervals are $[v_i,v_{i+5})$ for $1\le i\le 8$; they are 
the images under $T^{(i)}_{ab}$ of 
the interval $[(0,-1),v_5)$, $1\le i\le 8$.    These intervals are
entirely transient; that is, any orbit contains at most one point from each
of the eight intervals.

We consider  the interval $J= \RR^{+}[\bv, \bv') := \RR^{+}[\bv_9, \bv_5)$.
The break point $\bv^{*}= (0, -1).$
The first return map on $J_1= [\bv, \bv^{*})$ is $T_{ab}^{(5)}$
and has symbol sequence $(-+++-)$; its associated matrix
$M_1$ is
$$
M_1 = \left[ { {(a^3 -2 a)b^2 +(-2a^2 +2)b + a} \atop 
{(a^3 - 2a)b - a^2 +1}}~~
{{(-a^3 +2a)b +a^2 - 1} \atop {-a^3 + 2a}} \right] .
$$
The first return map on $J_2= [\bv^{*}, \bv')$ is $T_{ab}^{(9)}$
and has symbol sequence $(++++-++++)$; its associated matrix
$M_2$ is
$$
M_2= 
\frac{a^2-a-1}{a-1}M_1 + \frac{a^2(a-2)(a^2-2)^2}{(a-1)(a^2-a-1)}I.
$$
The argument of Theorem~\ref{th54} applies to show these two matrices
commute.
One checks that $|Tr(M_1)| > 2$ for $1 < a < \sqrt{2}$, hence the
associated matrices are hyperbolic.

We show that the rotation number
$$ 
r(T_{ab})=\frac{1}{5}.
$$
This follows because the
subsector $J'= \RR^{+}[(-1, 0), (-1,-1))$ is mapped into itself
by $T_J$ and therefore contains an attracting periodic point;
since it is in $J_1$ 
 this point has periodic symbol sequence $(-+++-)$ and therefore
has rotation number $1/5$. 

Lemma~\ref{le51} applies to show that on the sector $J$
the iterates of a given point lie on arcs of a hyperbola,
or on asymptotes.
One can show that the interval $J_1$ contains two asymptotes
of a hyperbola. Each of these defines a periodic point of order
$5$, one attracting and one repelling, defining the two periodic
orbits of period $5$ of $S_{ab}$. It follows that the 
first return map orbits
of all points diverge in norm to $+\infty$ under forward iterations,
except for those on the repelling fixed ray, which go to $+\infty$
under backwards iteration.

The points in $J$ not on an asymptote under iteration of
the first return map will
approach an asymptote, which is an attracting fixed point
of the first return map. For the map $T_{ab}$ it is an attracting
periodic orbit, and there are five copies of the asymptotes.
The orbit of any point not on a periodic orbit of $S_{ab}$ is
contained in a finite union of arcs of hyperbolas, even though
there are no invariant circles. Theorem~\ref{th55} does not apply to
bound the number of pieces, because the first return
map on  $J$
has three pieces; however the number of pieces appears to
be at most $23$. ~~~$\bsq$



\begin{thebibliography}{99}

\bibitem{AA68}
V. I. Arnold and A. Avez,
{\em Ergodic Problems of Classical Mechanics,}
W. A. Benjamin: New York 1968.

\bibitem{Ba75}
A. Baker,
{\em Transcendental Number Theory},
Cambridge Univ. Press, Cambridge 1975.

\bibitem{BBR95}
A. F. Beardon, S. R. Bullett and P. J. Rippon,
Periodic orbits of difference equations,
Proc. Roy. Soc. Edinburgh {\bf 125A} 91995), 657--674.

\bibitem{Fr68}
C. Froeschl\`{e},
\'{E}tude num\'{e}rique de transformations ponctuelles planes 
conservant les aires, C. R. Acad. Sci. Paris {\bf 266} (1968),
846--848.

\bibitem{He79}
M. Herman, 
Sur la conjugasion differentiable des diff\'{e}omorphismes du cercle.
Publ. Math. IHES {\bf 49 } (1979), 5--234.

\bibitem{He86}
M. Herman,
Sur les Courbes Invariantes par les Diff\'{e}omorphismes de l'Anneau. Vol. 2,
Ast\'{e}risque {\bf 144}, Soc. Math. de France: Paris 1986.


\bibitem{LR02aa}
J. C. Lagarias and E. Rains,
Dynamics of a  family of piecewise-linear area-preserving plane maps
 I. Rational rotation numbers, eprint: {\tt arxiv:math.DS/0301294}


\bibitem{LR02b}
J. C. Lagarias and E. Rains,
Dynamics of a  family of piecewise-linear area-preserving plane maps
 III. Cantor set spectra, eprint: {\tt arxiv:math.DS/0505103}

\bibitem{La93}
V. F. Lazutkin,
{\em KAM Theory and Semiclassical Approximations to Eigenfunctions},
Springer-Verlag: Berlin 1993.






\end{thebibliography}
\end{document}